\documentclass[reqno]{amsart}
\usepackage{amsrefs}
\usepackage{amssymb}
\usepackage{graphicx}
\usepackage{enumerate}
\usepackage{pst-tree}
\usepackage{pstricks-add}
\usepackage{url}

\usepackage{xcolor}
\newcommand{\upd}[1]{#1}

\newtheorem{thm}{Theorem}[section]
\newtheorem{prop}{Proposition}[section]

\numberwithin{equation}{section}

\DeclareMathOperator{\sech}{sech}
\newcommand{\abs}[1]{\left\lvert#1\right\rvert}
\newcommand{\ceiling}[1]{\left\lceil#1\right\rceil}
\newcommand{\floor}[1]{\left\lfloor#1\right\rfloor}
\newcommand{\symres}[2]{\langle#1\rangle_{#2}}

\newcommand{\tablerowsep}{2.5pt}
\newcommand\TS{\rule{0pt}{2.6ex}} 
\newcommand\BS{\rule[-1.2ex]{0pt}{0pt}} 

\allowdisplaybreaks

\begin{document}

\title{Wolstenholme and Vandiver primes}

\author[A.~R. Booker]{Andrew R. Booker}
\address{School of Mathematics, University of Bristol, Bristol, UK}
\email{andrew.booker@bristol.ac.uk}
\author[S. Hathi]{Shehzad Hathi}
\address{School of Science, The University of New South Wales, Canberra, Australia}
\email{s.hathi@student.adfa.edu.au}
\author[M.~J. Mossinghoff]{Michael J. Mossinghoff}
\address{Center for Communications Research, Princeton, NJ, USA}
\email{m.mossinghoff@idaccr.org}
\author[T.~S. Trudgian]{Timothy S. Trudgian}
\address{School of Science, The University of New South Wales, Canberra, Australia}
\thanks{Supported by Australian Research Council Future Fellowship FT160100094.}
\email{t.trudgian@adfa.edu.au}

\date\today

\begin{abstract}
A prime $p$ is a \textit{Wolstenholme prime} if $\binom{2p}{p}\equiv2$ mod $p^4$, or, equivalently, if $p$ divides the numerator of the Bernoulli number $B_{p-3}$; a \textit{Vandiver prime} $p$ is one that divides the Euler number $E_{p-3}$.
Only two Wolstenholme primes and eight Vandiver primes are known.
We increase the search range in the first case by a factor of \upd{ten}, and show that no additional Wolstenholme primes exist up to $\upd{10^{11}}$, and in the second case by a factor of \upd{twenty}, proving that no additional Vandiver primes occur up to \upd{this same bound}.
To facilitate this, we develop a number of new congruences for Bernoulli and Euler numbers mod $p$ that are favorable for computation, and we implement some highly parallel searches using GPUs.
\end{abstract}

\maketitle

\section{Introduction}\label{sectionIntroduction}

In 1862, Wolstenholme \cite{Wolstenholme} established that the following three congruences hold for every prime $p\geq5$:
\[
\sum_{0<k<p} k^{-2} \equiv 0 \mod p, \;\;
\sum_{0<k<p} k^{-1} \equiv 0 \mod p^2, \;\;
\binom{2p-1}{p-1} \equiv 1 \mod p^3.
\]
These results are closely related: if the stronger congruence obtained by replacing the modulus $p^m$ by $p^{m+1}$ holds in one of the expressions above for a particular prime $p$, then a similar strengthening for that same prime $p$ holds in all of them.
See \cite{Gardiner} for a proof, \cite{AebiCairns} for remarks on higher-order congruences, and \cite{MestrovicSurvey} for a survey on other generalizations and extensions of Wolstenholme's theorem.
A \textit{Wolstenholme prime} (first defined by McIntosh \cite{McIntosh}) is an odd prime $p$ for which any of the following congruences holds:
\begin{equation}\label{eqnWolst}
\sum_{0<k<p} k^{-2} \equiv 0 \mod p^2, \;\;
\sum_{0<k<p} k^{-1} \equiv 0 \mod p^3, \;\;
\binom{2p-1}{p-1} \equiv 1 \mod p^4.
\end{equation}
Many conditions equivalent to these are also well known.
For example, in the last congruence, the constant $2$ may be replaced by any integer $h$ not divisible by $p$.
It is also straightforward to rewrite that generalized congruence as
\[
\binom{hp}{p} \equiv h \mod p^4,
\]
which some authors employ.

Recall the \textit{Bernoulli numbers} $B_k$ may be defined by their exponential generating function,
\[
\frac{z}{e^z-1} = \sum_{k\geq0} B_k \frac{z^k}{k!}.
\]
These are rational numbers, and $B_{2k+1}=0$ for $k\geq1$.
 An important connection between certain binomial coefficients and Bernoulli numbers was established by Glaisher in 1900 (combining \cite[\S39]{Glaisher1} and \cite[\S8]{Glaisher2}), who proved that
\begin{equation}\label{eqnGlaisher}
\binom{hp-1}{p-1} \equiv 1 - \frac{h(h-1)}{3} p^3 B_{p-3} \mod p^4
\end{equation}
for any prime $p\geq5$ and any positive integer $h$.
Thus, $p$ is a Wolstenholme prime if and only if $p$ divides the numerator of the Bernoulli number $B_{p-3}$.
This condition provides the first observation toward a more computationally attractive test for Wolstenholme primes, since it requires computation of a residue mod $p$, rather than a higher power of $p$.

Glaisher's observation implies that a Wolstenholme prime is a particular type of \textit{irregular prime}.
Recall that a prime $p$ is said to be \textit{regular} if $p$ does not divide the class number of the cyclotomic field $\mathbb{Q}(\zeta_p)$, where $\zeta_p=e^{2\pi i/p}$.
Equivalently, $p$ is regular if $p$ does not divide the numerator of any of the Bernoulli numbers $B_{2k}$ with $0<2k\leq p-3$.
Regular primes were first defined in the context of work on Fermat's Last Theorem, as Kummer proved in 1847 that no nontrivial solutions with exponent $p$ exist in this famous problem if $p$ is a regular prime.
The smallest irregular prime is $37$, and while it is known that there exist infinitely many irregular primes (see for instance \citelist{\cite{Carlitz54}\cite{LPMP15}}), it remains unknown if the number of regular primes is infinite.
We refer the reader to \cites{Edwards,Washington} for more information on regular primes and their connections to various topics in number theory.

The \textit{Euler numbers} $E_k$ may also be defined by their exponential generating function,\footnote{Some authors instead use the series for the \textsl{hyperbolic} secant, $\sech z = \frac{2e^z}{e^{2z}+1}$, to define the Euler numbers; this alters the sign of the terms with index congruent to $2$ mod $4$. Results cited from the literature that employ that formulation are suitably translated here to comport with the definition we employ.}
\[
\sec z = \sum_{k\geq0} E_k \frac{z^k}{k!}.
\]
Like the Bernoulli numbers, the Euler numbers with odd index are $0$ (including $E_1$ this time), but unlike the $B_k$, the Euler numbers are all integers.
A prime $p$ is said to be \textit{$E$-regular} if $p$ does not divide any of the Euler numbers $E_{2k}$ with $0<2k\leq p-3$.
The smallest $E$-irregular prime is $19$, as $19\mid E_{10}$, and while it is known that there exist infinitely many $E$-irregular primes \citelist{\cite{Carlitz54}\cite{LPMP15}}, it is also not known if the number of $E$-regular primes is infinite.
$E$-irregular primes first arose in the study of Fermat's Last Theorem as well, and they are pertinent in other problems in number theory too, including the topic of class numbers of cyclotomic fields.
For example, Gut \cite{Gut} proved that $p$ does not divide the class number of the cyclotomic field $\mathbb{Q}(\zeta_{4p})$ if and only if $B_{2k}\not\equiv0$ and $E_{2k}\not\equiv0$ mod $p$ for $0<2k\leq p-3$.

The extremal case $E_{p-3}\equiv0$ mod $p$ has witnessed particular interest.
In 1940, Vandiver \cite{Vandiver40} showed that any solution with exponent $p$ in the first case of Fermat's Last Theorem (so a nontrivial solution in positive integers to $x^p+y^p=z^p$ where $p\nmid xyz$) necessarily has $p\mid E_{p-3}$.
More recently, Z.-H. Sun \cite{ZHSun08} and Z.-W. Sun \cite{ZWSun11} established a number of congruences regarding Bernoulli and Euler numbers, in many of which the residues of $B_{p-3}$ and $E_{p-3}$ mod $p$ play an important role, and several of these congruences become stronger in the case of a Wolstenholme or Vandiver prime. 

If $p$ is irregular, and $B_{2r}\equiv0$ mod $p$ for an even integer $2r\leq p-3$, then $(p,2r)$ is known as an \textit{irregular pair}.
A Wolstenholme prime is thus the extremal case where $(p,p-3)$ forms an irregular pair.
Only two Wolstenholme primes are known, and both were discovered during broader searches for irregular primes.
The smallest Wolstenholme prime is $p_1=16843$.
This seems to have been first discovered by Selfridge and Pollack, who spoke on their work searching for irregular primes up to $25000$ at a meeting of the American Mathematical Society in 1964 \cite{SelfridgePollack}.
While details from their search are not readily available, in 1975 Johnson \cite{Johnson75} reported on an extension of the work of Selfridge and Pollack, verifying their work, and noted in particular that $p_1$ was the only prime less than $30000$ where $(p,p-3)$ forms an irregular pair.

The second Wolstenholme prime is $p_2=2124679$, and was discovered in 1993 by Buhler et al.\ \cite{BCEM93} during a search for irregular primes $p<4\cdot10^6$.
In 1995 McIntosh \cite{McIntosh} reported finding no additional examples with $p<2\cdot10^8$, after more extensive searches that targeted Wolstenholme primes in particular, and then (with Roettger) extended this search in 2007 \cite{McIntoshRoettger}  to $p<10^9$.
More recently, McIntosh continued this search to cover $p<10^{10}$ \cite{McIntoshFQ} and found no additional examples.
Broader searches for irregular primes have continued as well, culminating most recently in the work of Hart, Harvey, and Ong \cite{HHO} of 2017, who determined all irregular primes $p<2^{31}$, along with their index of irregularity, that is, the number of irregular pairs attached to each such prime.
They found for instance that the maximum index of irregularity over this range is $9$, attained by the single prime $p=1767218027$.

We likewise say $(p,2r)$ is an \textit{$E$-irregular pair} if $E_{2r}\equiv0$ mod $p$, so that a Vandiver prime has the property that $(p,p-3)$ forms an $E$-irregular pair.
It appears that the last census of $E$-irregular primes and indices of $E$-irregularity was performed in 1978 by Ernvall and Mets\"{a}nkyl\"{a} \citelist{\cite{Ernvall79}\cite{EM78}}, who checked primes $p<10^4$, motivated by questions regarding Iwasawa invariants of the cyclotomic field $\mathbb{Q}(\zeta_{4p})$.
They found for example that the maximum index of $E$-irregularity over this range was $5$, attained by $p=5783$ alone.

A number of Vandiver primes were known prior to this work.
The first two, $p=149$ and $241$, were first noted by Ernvall and Mets\"{a}nkyl\"{a} \cite{EM78} (and later rediscovered by Z.-W. Sun \cite{ZWSun11}).
Three more Vandiver primes, $p=2124679$, $16467631$, and $17613227$, were found by Cosgrave and Dilcher in 2013 \cite{CosgraveDilcher} (who called these integers \textit{$E$-primes}), using a congruence from \cite{ELehmer} to search up to $5\cdot10^7$.
They showed that these primes arise in a natural way in a problem regarding modular sums of certain reciprocals, proving that
\[
\sum_{\substack{0<k\leq\floor{n/4}\\\gcd(k,n)=1}} k^{-2} \equiv 0 \mod n
\]
precisely when $n=45$ or $n$ is divisible by a Vandiver prime.
Cosgrave and Dilcher also reported that in 2012 McIntosh found three additional such primes in a search to $3\cdot10^9$: $p=327784727$, $426369739$, and $1062232319$.
More recently, McIntosh extended this to $5\cdot10^9$ \cite{McIntoshFQ}, finding no additional examples.
We also report that in 2014 Me\v{s}trovi\'c \cite{Mestrovic14} employed a congruence of Z.-H. Sun from \cite{ZHSun08} involving the computation of certain harmonic numbers modulo $p^2$ to compute the Vandiver primes $p<10^7$.

In the absence of any special structure, one might expect the values of $B_{p-3}$ or $E_{p-3}$ to be uniformly distributed modulo $p$.\footnote{We remark that some structure does appear in a similar problem.
The Ankeny--Artin--Chowla conjecture asserts that $p\nmid B_{(p-1)/2}$ for every prime $p\equiv1$ mod $4$.
While this problem is open in general (and has been verified computationally for $p<2\cdot10^{11}$ \cite{VTW}), this conjecture is known to hold for primes of the form $p=n^2+1$ or $n^2+4$ \cite{Agoh}.}
Numerical evidence weighs in support of this: see Figures~\ref{figBernHist} and~\ref{figEulerHist} later in this paper.
If this is true, then as McIntosh \cite{McIntosh} has pointed out, one would expect that infinitely many of each of these kinds of primes exist, since the sum of the reciprocals of the primes diverges.
In this case, the expected number of either type of prime up to $x$ would grow approximately as $\log\log x$.
One might define a \textit{super-Wolstenholme} prime by requiring a congruence in \eqref{eqnWolst} to hold modulo a higher power of $p$, and similarly for a \textit{super-Vandiver} prime.
Similar heuristic considerations lead one to expect that only a finite number of such primes exist.
In \cite{McIntosh} McIntosh conjectured that no super-Wolstenholme primes exist; this seems likely to be true as well for super-Vandiver primes.

In this article, we describe searches for Wolstenholme and Vandiver primes that reach eight times farther than any prior search in the former case, and twelve times farther in the latter.
We establish that no additional examples of either type occur across the ranges we check.

\begin{thm}\label{thmWVknown}
The Wolstenholme primes less than \upd{$10^{11}$} are
$16843$ and $2124679$;
the Vandiver primes less than \upd{$10^{11}$} are
$149$, $241$, $2124679$, $16467631$, $17613227$, $327784727$, $426369739$, and $1062232319$.
\end{thm}

There are two principal components to our method which allowed us to extend these searches by substantial factors.
First, we obtain new congruences for Bernoulli and Euler numbers modulo a prime $p$, which allow us to compute residues of these numbers with greater efficiency.
See \eqref{eqnB6}, \eqref{eqnB9}, \eqref{eqnB16}, \eqref{eqnB22}, and \eqref{eqnB30} for successively more efficient (and more complicated) congruences for the Bernoulli numbers, and similarly see \eqref{eqnE3}, \eqref{eqnE5}, \eqref{eqnE9}, \eqref{eqnE16}, \eqref{eqnE24}, and \eqref{eqnE33} for our congruences for the Euler numbers.
For the Bernoulli numbers, our congruences extend and improve on some results of Tanner and Wagstaff from 1987 \cite{TannerWagstaff}.
Second, we employ graphics processing units (GPUs) to implement a heavily parallelized approach.

\upd{We also answer in the affirmative a question raised by Tanner and Wagstaff regarding the existence of congruences for Bernoulli numbers mod $p$ involving fewer than $\epsilon p$ terms, for arbitrarily small $\epsilon>0$---see Theorem~\ref{thmSmallCost}.
This result applies equally well to the computation of Euler numbers mod $p$.}

This article is organized in the following way.
Section~\ref{sectionBernoulli} derives a number of new congruences for Bernoulli numbers which are advantageous for computing, particularly for large primes.
Section~\ref{sectionComputationsW} then describes our search for Wolstenholme primes using these congruences and summarizes the results of our searches in this problem up to \upd{$10^{11}$}.
Section~\ref{sectionEuler} establishes some new congruences amenable for our computations on Euler numbers, and Section~\ref{sectionComputationsV} describes our search for Vandiver primes and summarizes our results up to \upd{this same bound}.
\upd{Section~\ref{sectionSmallCost} establishes that congruences with fewer than $\epsilon p$ terms exist for computing residues of Bernoulli or Euler numbers, for arbitrary $\epsilon>0$, and} Section~\ref{sectionBounds} adds some observations regarding derivations of favorable congruences.

\section{Congruences for Bernoulli numbers}\label{sectionBernoulli}

The Bernoulli numbers are known to satisfy a number of congruences involving sums of powers of consecutive integers, beginning with Voronoi's well-known congruence of 1889 (see \cite[Chap.\ IX]{UspenskyHeaslet}).
Some of these congruences have been particularly instrumental in calculations involving Bernoulli numbers.
We review some of the more important relations.
For real numbers $x<y$ in $[0,1]$, an integer $\ell$, and a prime $p$, let
\begin{equation}\label{eqnDefnS}
S_\ell(x,y) = S_{\ell,p}(x,y) = \sum_{xp<s<yp} s^\ell \mod p,
\end{equation}
and, for positive integers $a$, $b$, $c$ and a positive integer $k$, let
\[
C_k(a,b,c) = C_{k,p}(a,b,c) = \frac{a^{p-2k} + b^{p-2k} - c^{p-2k} - 1}{4k}.
\]
With this notation, in 1930 Stafford and Vandiver \cite[eq.\ (9)]{StaffordVandiver}] (see also \cite[Thm.\ 1]{TannerWagstaff}) established that
\begin{equation}\label{eqnSV}
\textstyle
C_k(3,4,6) B_{2k} \equiv S_{2k-1}(\frac{1}{6}, \frac{1}{4}) \mod p 
\end{equation}
whenever $p\geq5$ and $p-1 \nmid 2k$.
Note that this allows computing the residue of $B_{2k}$ mod $p$ using a sum with approximately $p/12$ terms.
Also, in 1937, Vandiver \cite[eq.\ (18)]{Vandiver37} proved that
\begin{equation}\label{eqnVandiver}
\textstyle
C_k(2,5,6) B_{2k} \equiv S_{2k-1}(\frac{1}{6}, \frac{1}{5}) + S_{2k-1}(\frac{1}{3}, \frac{2}{5}) \mod p  
\end{equation}
when $p\geq7$ and $p-1 \nmid 2k$.
For convenience, we define the \textit{cost} of a congruence such as \eqref{eqnSV} or \eqref{eqnVandiver}, where in general the right side consists of a linear combination of sums $S_\ell(x_1,y_1)$, \ldots, $S_\ell(x_m,y_m)$, as
\[
p\sum_{i=1}^m(y_i-x_i),
\]
as this is the approximate number of terms in the sums comprising the right side of the congruence.
The cost of Vandiver's congruence \eqref{eqnVandiver} is thus $p(\frac{1}{5}-\frac{1}{6}+\frac{2}{5}-\frac{1}{3})=p/10$, and the cost of \eqref{eqnSV} is $p/12$.
However, it transpires that \eqref{eqnVandiver} is more important for our application.

In 1987, Tanner and Wagstaff \cite{TannerWagstaff} developed a number of new congruences for Bernoulli numbers.
They proved for example that
\begin{equation}\label{eqnTW1}
C_k(2,b,b+1) B_{2k} \equiv \sum_{m=1}^{\floor{b/2}} S_{2k-1}\left(\frac{m}{b+1},\frac{m}{b}\right) \mod p
\end{equation}
whenever $p$ is prime, $b\geq2$, $p>b+1$, $k\geq1$, and $p-1\nmid 2k$.
This rule has cost $\frac{\floor{b/2}(\floor{b/2}+1)}{2b(b+1)}p$, so asymptotically $p/8$, and we note that the case $b=5$ recovers \eqref{eqnVandiver}.
Tanner and Wagstaff also noted three simple transformations that preserve the value of $S_\ell(x,y)$ mod $p$, which we summarize here, in a slightly generalized form.

\begin{prop}[Tanner and Wagstaff \cite{TannerWagstaff}]\label{propTW}
The following three properties hold for $S_\ell(x,y)$, defined in \eqref{eqnDefnS}.
\begin{enumerate}[(a)]
\item Separation: If $0\leq x<y<z\leq 1$ and $yp$ is not an integer, then
\[
S_\ell(x,z) = S_\ell(x,y) + S_\ell(y,z).
\]
\item Reflection: $S_\ell(x,y) \equiv (-1)^\ell S_\ell(1-y,1-x) \mod p$.
\item Subdivision: If $d$ is a positive integer and $p\nmid d$, then
\[
S_\ell(x,y) \equiv d^\ell \sum_{i=0}^{d-1} S_\ell\left(\frac{x+i}{d}, \frac{y+i}{d}\right) \mod p.
\]
\end{enumerate}
\end{prop}

We include the short proof here for the convenience of the reader.

\begin{proof}
The first two properties are immediate from the definition of $S_\ell(x,y)$.
For the third, since $p\nmid d$, every integer $s$ can be expressed in a unique way as $s=qd-ip$ for some integer $i\in[0,d)$.
Therefore,
\begin{align*}
S_\ell(x,y) &= \sum_{i=0}^{d-1} \sum_{xp < qd-ip < yp} (qd-ip)^\ell\\
&\equiv \sum_{i=0}^{d-1} \sum_{xp < qd-ip < yp} (qd)^\ell \mod p\\
&\equiv d^\ell \sum_{i=0}^{d-1} S_\ell\left(\frac{x+i}{d}, \frac{y+i}{d}\right) \mod p.\qedhere
\end{align*}
\end{proof}

One may employ the transformations of Proposition~\ref{propTW} to create new congruences for Bernoulli numbers from existing ones.
While each transformation preserves the cost of the sums involved, one might hope that subdivision and reflection might create new expressions with some overlapping sums, so that separation may be employed to collect like terms, and thus produce congruences with reduced cost, at the price of some more complicated coefficients.
For example, as noted in \cite{TannerWagstaff}, if one applies subdivision with $d=2$ to the second term of \eqref{eqnVandiver}, and then reflection to the resulting term having $x>1/2$, then one obtains
\begin{equation}\label{eqnB2}
\begin{split}
\textstyle C_k(2,5,6) B_{2k}
&\textstyle \equiv S_{2k-1}(\frac{1}{6}, \frac{1}{5}) + 2^{2k-1}\bigl(S_{2k-1}(\frac{1}{6}, \frac{1}{5}) + S_{2k-1}(\frac{2}{3}, \frac{7}{10})\bigr)\\
&\textstyle \equiv (1+2^{2k-1})S_{2k-1}(\frac{1}{6}, \frac{1}{5}) - 2^{2k-1}S_{2k-1}(\frac{3}{10}, \frac{1}{3}) \mod p,
\end{split}
\end{equation}
provided $p\geq7$ and $p-1 \nmid 2k$.
This new expression has two sums, each now possessing modestly more complicated coefficients, but the cost of this congruence has dropped to $p/15$.
Tanner and Wagstaff also derived a more complicated congruence for Bernoulli numbers with more sums and lower cost.
By starting with \eqref{eqnTW1} in the case $b=9$, and employing five prescribed subdivision steps in sequence, using reflection and separation as appropriate for simplification, they found a congruence for $B_{2k}$ mod $p$ having seven sums with cost $p/18$.
They also reported that beginning with the congruence \eqref{eqnB2} and applying two particular subdivision operations, together with appropriate simplifications, produces an expression with cost $p/19.2$.
While the resulting congruence is not reported, we find that it has eight sums.
Finally, it is mentioned in the same article that the authors could construct a more complicated congruence with approximate cost $p/22$, but its construction is not specified.

We use Proposition~\ref{propTW} to search for congruences for Bernoulli numbers mod $p$ that are particularly favorable for computations.
Certainly, it is advantageous to employ a congruence with small cost, but one must also be concerned with the number of sums appearing in the expressions, as well as the size and complexity of their coefficients, as each of these contributes to overhead costs in computations.
It may also be helpful for each sum in a congruence to have a comparable number of elements, for balancing loads in parallel environments.
In addition, for larger primes we may be able to tolerate greater overhead costs, so we therefore seek a family of congruences for Bernoulli numbers mod $p$: for a range of small positive integers $m$, we would like an expression for $B_{2k}$ having exactly $m$ sums with cost $p/r_m$, with $r_m$ as large as possible.

The case $m=1$ was the subject of a search by Wagstaff \cite{Wagstaff}, who found nothing better than the value $r_1=12$ exhibited by \eqref{eqnSV}, as well as by the well-known relation (see, e.g., \cite[ex.\ IX.6]{UspenskyHeaslet})
\begin{equation}\label{eqnV12}
\textstyle
C_k(2,3,4) B_{2k} \equiv S_{2k-1}(\frac{1}{4}, \frac{1}{3}) \mod p,
\end{equation}
valid for primes $p\geq5$ when $p-1 \nmid 2k$.
This latter relation follows from Voronoi's congruence, and in fact is the case $b=3$ in \eqref{eqnTW1}.
We searched for improved values for some larger $m$ by using Proposition~\ref{propTW}, with a number of starting relations.

We employed two strategies in searching for useful congruences: exhaustive searches, and heuristic searches.
We describe each of these below, along with its principal results.

\subsection{Exhaustive searches}\label{subsecExhaustive}

In an exhaustive search, we supply a starting congruence such as \eqref{eqnSV}, along with a list of positive integers $d_1, \ldots, d_\ell$, and we construct all possible congruences that can be formed with a sequence of at most $\ell$ subdivision operations starting from our initial congruence, where the $i$th subdivision employs a value of $d$ satisfying $2\leq d\leq d_i$.
We make a number of runs with this method using different initial lists $d_1,\ldots,d_\ell$, some with a small value of $\ell$ like $\ell=4$ and larger values for the $d_i$ such as $6$ or $8$, and some with a larger value of $\ell$ like $\ell=7$, where we limit the $d_i$ to smaller values such as $2$ or $3$.
After each subdivision, we use the reflection and separation properties to simplify our expressions as much as possible.
Over several runs with different parameters we maintain a list that stores for each $m$ the smallest cost of a congruence having $m$ sums that was constructed from the given starting relation, along with information on how to construct it.
For starting relations, we employed each of the congruences \eqref{eqnSV}, \eqref{eqnVandiver}, \eqref{eqnTW1} with $b\leq13$, \eqref{eqnV12}, as well as the well-known relation
\begin{equation}\label{eqnV10}
\textstyle
C_k(4,5,8) B_{2k} \equiv
S_{2k-1}(\frac{1}{8}, \frac{1}{5}) + S_{2k-1}(\frac{3}{8}, \frac{2}{5})\mod p,
\end{equation}
which is valid for $p\geq7$ and $p-1 \nmid 2k$ and has cost $p/10$.
This last expression also follows from Voronoi's congruence \cite[ex.\ IX.6]{UspenskyHeaslet}.
In addition, in order to search more effectively at larger values of $m$, we would start this procedure at a favorable congruence observed at a smaller value of $m$.
For example, many improved values for $m\geq18$ were found by starting a search at the best known congruence at $m=16$.

Table~\ref{tblBestCost} summarizes the results of our searches using this strategy, displaying the minimal cost discovered among congruences having exactly $m$ sums, for a range of values of $m$.
In each case, the best example we found was constructed using Vandiver's congruence \eqref{eqnVandiver} as its starting point.
We therefore assume $p\geq7$ and $p-1 \nmid 2k$ in all congruences displayed in the remainder of this section.
We provide details only on the congruences from Table~\ref{tblBestCost} that were employed in our searches for Wolstenholme primes, described in Section~\ref{sectionComputationsW}.

\begin{table}[tb]
\caption{For $m\leq24$, the smallest known cost for a congruence for Bernoulli numbers mod $p$ having exactly $m$ sums is $p/r_m$.
For $m>1$, each was constructed using the method of Section~\ref{subsecExhaustive} using \eqref{eqnVandiver} as the initial congruence.}\label{tblBestCost}
\begin{tabular}{|c|c||c|c||c|c|}\hline
$m$ & $r_m$ & $m$ & $r_m$ & $m$ & $r_m$\\[\tablerowsep]\hline
\TS$1$ & $12$ & $9$ & $20$ & $17$ & $\frac{1080}{43} = 25.11\ldots$\\[\tablerowsep]
\TS$2$ & $15$ & $10$ & $\frac{1200}{59} = 20.33\ldots$ & $18$ & $\frac{180}{7} = 25.71\ldots$\\[\tablerowsep]
\TS$3$ & $15$ & $11$ & $\frac{360}{17} = 21.17\ldots$ & $19$ & $\frac{180}{7} = 25.71\ldots$\\[\tablerowsep]
\TS$4$ & $\frac{120}{7} = 17.14\ldots$ & $12$ & $\frac{108}{5} = 21.6$ & $20$ & $\frac{810}{31} = 26.12\ldots$\\[\tablerowsep]
\TS$5$ & $\frac{240}{13} = 18.46\ldots$ & $13$ & $\frac{45}{2} = 22.5$ & $21$ & $\frac{1080}{41} = 26.34\ldots$\\[\tablerowsep]
\TS$6$ & $\frac{96}{5} = 19.2$ & $14$ & $\frac{1080}{47} = 22.97\ldots$ & $22$ & $\frac{80}{3} = 26.66\ldots$\\[\tablerowsep]
\TS$7$ & $\frac{960}{49} = 19.59\ldots$ & $15$ & $\frac{540}{23} = 23.47\ldots$ & $23$ & $\frac{3240}{121} = 26.77\ldots$\\[\tablerowsep]
\TS$8$ & $\frac{1920}{97} = 19.79\ldots$ & $16$ & $24$ & $24$ & $\frac{1620}{59} = 27.45\ldots$\\[\tablerowsep]\hline
\end{tabular}
\end{table}

For $m=2$, the best example we find is in fact \eqref{eqnB2}, where $r_2=15$.

For $m=6$, we find a congruence with cost $p/19.2$, so this requires two fewer sums than the congruence with the same cost reported in \cite{TannerWagstaff}.
To obtain this, starting with \eqref{eqnB2} one subdivides the $S_{2k-1}(3/10,1/3)$ term with $d=2$, simplifies appropriately, then subdivides $S_{2k-1}(1/3, 7/20)$ using $d=2$, then the $S_{2k-1}(13/40, 1/3)$ term with $d=2$, and finally subdivides the resulting sum $S_{2k-1}(1/3, 27/80)$ with $d=2$.
Writing $t=2k-1$ for economy, this produces the six-term congruence
\begin{equation}\label{eqnB6}
\begin{split}
\textstyle
C_k(2,{}&5,6)
\textstyle
B_{2k} \equiv 
- 2^{2t} S_t(\frac{3}{20}, \frac{13}{80})
- (2^{2t} + 2^{4t}) S_t(\frac{13}{80}, \frac{1}{6})\\
&\textstyle
+ (1 + 2^t + 2^{3t} + 2^{5t}) S_t(\frac{1}{6}, \frac{27}{160})
+ (1 + 2^t + 2^{3t}) S_t(\frac{27}{160}, \frac{7}{40})\\
&\textstyle
+ (1 + 2^t) S_t(\frac{7}{40}, \frac{1}{5})
- 2^{5t} S_t(\frac{53}{160}, \frac{1}{3})
\mod p.
\end{split}
\end{equation}
Here, the cost of each sum is respectively $p/80$, $p/240$, $p/480$, $p/160$, $p/40$, and $p/480$, so this example has $r_6=96/5$.

We remark that intermediate expressions generated when creating \eqref{eqnB6} in this way produce two more of the entries in Table~\ref{tblBestCost}.
Applying just the first two subdivisions described above to \eqref{eqnB2} produces our congruence with $m=4$ sums and cost $7p/120$; applying the third subdivision to this manufactures our example at $m=5$ with cost $13p/240$.

We also find a particularly well-balanced example with small cost at $m=9$.
Again starting with Vandiver's congruence \eqref{eqnVandiver}, we first subdivide $S_t(1/6,1/5)$ with $d=3$, then subdivide the resulting sums $S_t(4/15, 5/18)$ and $S_t(7/18, 2/5)$ using $d=2$, and finally subdivide $S_t(11/36, 1/3)$ with $d=5$.
This produces the nine-term congruence
\begin{equation}\label{eqnB9}
\begin{split}
 C_k(2,{}&5,6)
 \textstyle
 B_{2k} \equiv 
 (3^t + 6^t) S_t(\frac{1}{18}, \frac{11}{180})
+ (3^t + 6^t - 10^t) S_t(\frac{11}{180}, \frac{1}{15})\\
&\textstyle - (6^t - 10^t + 12^t) S_t(\frac{2}{15}, \frac{5}{36})
+ (6^t + 12^t) S_t(\frac{7}{36}, \frac{1}{5})\\
&\textstyle - 10^t S_t(\frac{47}{180}, \frac{4}{15})
- (2^t + 6^t + 12^t) S_t(\frac{3}{10}, \frac{11}{36})
 + 10^t S_t(\frac{1}{3}, \frac{61}{180})\\
&\textstyle+ (6^t + 12^t) S_t(\frac{13}{36}, \frac{11}{30})
 - 10^t S_t(\frac{83}{180}, \frac{7}{15})
\mod p,
\end{split}
\end{equation}
where again we write $t$ for $2k-1$.
Each sum here has the same cost, $p/180$, so this expression achieves $r_9=20$.
Figure~\ref{m9congruencegraph} in Section~\ref{sectionBounds} depicts the transformation of \eqref{eqnVandiver} into \eqref{eqnB9}.

We provide details on two additional congruences from Table~\ref{tblBestCost}, which are employed in the searches described in the next section.
For an expression with $m=16$ sums and cost $p/24$, starting from \eqref{eqnB9} we subdivide both the $S_t(7/36, 1/5)$ and the $S_t(3/10, 11/36)$ sums using $d=3$, then subdivide $S_t(47/180, 4/15)$ with $d=2$, and finally subdivide the resulting sums $S_t(47/360, 2/15)$ and $S_t(83/180, 7/15)$ with $d=2$.
This yields
\begin{equation}\label{eqnB16}
\begin{split}
C_k(2,{}&5,6)
\textstyle
B_{2k} \equiv 
(3^t + 6^t) S_t(\frac{1}{18}, \frac{11}{180})
+ (3^t + 6^t - 10^t) S_t(\frac{11}{180}, \frac{7}{108})\\
&\textstyle
+ (3^t + 6^t - 10^t + 18^t + 36^t) S_t(\frac{7}{108}, \frac{47}{720})\\
&\textstyle
+ (3^t + 6^t - 10^t + 18^t + 36^t - 40^t) S_t(\frac{47}{720}, \frac{1}{15})\\
&\textstyle
- (6^t + 18^t + 36^t) S_t(\frac{1}{10}, \frac{11}{108})
- (6^t - 10^t + 12^t) S_t(\frac{2}{15}, \frac{5}{36})\\
&\textstyle
- 20^t S_t(\frac{83}{360}, \frac{25}{108})
+ (6^t + 18^t - 20^t + 36^t) S_t(\frac{25}{108}, \frac{7}{30})\\
&\textstyle
- (18^t - 20^t + 36^t) S_t(\frac{4}{15}, \frac{29}{108})
+ 20^t S_t(\frac{29}{108}, \frac{97}{360})\\
&\textstyle
+ 10^t S_t(\frac{1}{3}, \frac{61}{180})
+ (6^t + 12^t) S_t(\frac{13}{36}, \frac{11}{30})
+ 20^t S_t(\frac{11}{30}, \frac{133}{360})\\
&\textstyle
+ (18^t + 36^t) S_t(\frac{43}{108}, \frac{2}{5})
- (6^t + 18^t + 36^t - 40^t) S_t(\frac{13}{30}, \frac{313}{720})\\
&\textstyle
- (6^t + 18^t + 36^t) S_t(\frac{313}{720}, \frac{47}{108})
\mod p.
\end{split}
\end{equation}
Here, while the costs of the individual sums vary considerably, from $p/2160$ to $p/180$,
the cost of the entire rule is just $p/24$, so exactly half that of \eqref{eqnSV}.

We also record the derivation of the $22$-term congruence referenced in Table~\ref{tblBestCost}.
Starting from \eqref{eqnB16}, subdivide $S_t(1/3, 61/180)$ with $d=3$, and then subdivide each of $S_t(1/9, 61/540)$, $S_t(4/15, 29/108)$, $S_t(29/108, 97/360)$, and $S_t(479/1080, 4/9)$ in turn using $d=2$.
This yields
\begin{equation}\label{eqnB22}
\begin{split}
&\textstyle
C_k(2, 5, 6) B_{2k} \equiv 
(3^t + 6^t + 60^t) S_t(\frac{1}{18}, \frac{61}{1080})
+ (3^t + 6^t) S_t(\frac{61}{1080}, \frac{11}{180})\\
&\textstyle\;
+ (3^t + 6^t - 10^t) S_t(\frac{11}{180}, \frac{7}{108})
+ (3^t + 6^t - 10^t + 18^t + 36^t) S_t(\frac{7}{108}, \frac{47}{720})\\
&\textstyle\;
+ (3^t + 6^t - 10^t + 18^t + 36^t - 40^t) S_t(\frac{47}{720}, \frac{1}{15})\\
&\textstyle\;
- (6^t + 18^t + 36^t) S_t(\frac{1}{10}, \frac{11}{108})\\
&\textstyle\;
- (6^t - 10^t + 12^t + 36^t - 40^t + 72^t) S_t(\frac{2}{15}, \frac{29}{216})\\
&\textstyle\;
- (6^t - 10^t + 12^t - 40^t) S_t(\frac{29}{216}, \frac{97}{720})\\
&\textstyle\;
- (6^t - 10^t + 12^t) S_t(\frac{97}{720}, \frac{5}{36})
- 30^t S_t(\frac{119}{540}, \frac{479}{2160})
- (30^t + 120^t) S_t(\frac{479}{2160}, \frac{2}{9})\\
&\textstyle\;
- 20^t S_t(\frac{83}{360}, \frac{25}{108})
+ (6^t + 18^t - 20^t + 36^t) S_t(\frac{25}{108}, \frac{7}{30})
+ 120^t S_t(\frac{5}{18}, \frac{601}{2160})\\
&\textstyle\;
+ (6^t + 12^t) S_t(\frac{13}{36}, \frac{263}{720})
+ (6^t + 12^t - 40^t) S_t(\frac{263}{720}, \frac{79}{216})\\
&\textstyle\;
+ (6^t + 12^t + 36^t - 40^t + 72^t) S_t(\frac{79}{216}, \frac{11}{30})
+ 20^t S_t(\frac{11}{30}, \frac{133}{360})\\
&\textstyle\;
+ (18^t + 36^t) S_t(\frac{43}{108}, \frac{2}{5})
- (6^t + 18^t + 36^t - 40^t) S_t(\frac{13}{30}, \frac{313}{720})\\
&\textstyle\;
- (6^t + 18^t + 36^t) S_t(\frac{313}{720}, \frac{47}{108})
+ 30^t S_t(\frac{4}{9}, \frac{241}{540})
\mod p,
\end{split}
\end{equation}
which has cost $3p/80$.

Finally, we describe an additional congruence for Bernoulli numbers with $30$ sums that is employed in some of our searches: from \eqref{eqnB22}, subdivide $S_t(43/108, 2/5)$ using $d=6$, and then apply this transformation with $d=2$ to each of $S_t(4/15, 173/648)$, $S_t(4/9, 241/540)$, and $S_t(299/1080, 5/18)$ in turn.
This produces
\begin{equation}\label{eqnB30}
\begin{split}
&\scriptstyle
C_k(2, 5, 6) B_{2k} \equiv 
(3^t + 6^t + 60^t) S_t(\frac{1}{18}, \frac{61}{1080})
+ (3^t + 6^t) S_t(\frac{61}{1080}, \frac{11}{180})
+ (3^t + 6^t - 10^t) S_t(\frac{11}{180}, \frac{7}{108})\\
&\scriptstyle\;
+ (3^t + 6^t - 10^t + 18^t + 36^t) S_t(\frac{7}{108}, \frac{47}{720})
+ (3^t + 6^t - 10^t + 18^t + 36^t - 40^t) S_t(\frac{47}{720}, \frac{43}{648})\\
&\scriptstyle\;
+ (3^t + 6^t - 10^t + 18^t + 36^t - 40^t + 108^t + 216^t) S_t(\frac{43}{648}, \frac{1}{15})
- (6^t + 18^t + 36^t + 108^t + 216^t) S_t(\frac{1}{10}, \frac{65}{648})\\
&\scriptstyle\;
- (6^t + 18^t + 36^t) S_t(\frac{65}{648}, \frac{11}{108})
- (6^t - 10^t + 12^t + 36^t - 40^t + 72^t + 216^t + 432^t) S_t(\frac{2}{15}, \frac{173}{1296})\\
&\scriptstyle\;
- (6^t - 10^t + 12^t + 36^t - 40^t + 72^t) S_t(\frac{173}{1296}, \frac{29}{216})
- (6^t - 10^t + 12^t - 40^t) S_t(\frac{29}{216}, \frac{97}{720})\\
&\scriptstyle\;
- (6^t - 10^t + 12^t) S_t(\frac{97}{720}, \frac{299}{2160})
- (6^t - 10^t + 12^t + 120^t) S_t(\frac{299}{2160}, \frac{5}{36})
- 30^t S_t(\frac{119}{540}, \frac{479}{2160})\\
&\scriptstyle\;
- (30^t + 120^t) S_t(\frac{479}{2160}, \frac{2}{9})
+ 60^t S_t(\frac{2}{9}, \frac{241}{1080})
- 20^t S_t(\frac{83}{360}, \frac{25}{108})
+ (6^t + 18^t - 20^t + 36^t) S_t(\frac{25}{108}, \frac{151}{648})\\
&\scriptstyle\;
+ (6^t + 18^t - 20^t + 36^t + 108^t + 216^t) S_t(\frac{151}{648}, \frac{7}{30})
+ 120^t S_t(\frac{5}{18}, \frac{601}{2160})
+ (6^t + 12^t + 120^t) S_t(\frac{13}{36}, \frac{781}{2160})\\
&\scriptstyle\;
+ (6^t + 12^t) S_t(\frac{781}{2160}, \frac{263}{720})
+ (6^t + 12^t - 40^t) S_t(\frac{263}{720}, \frac{79}{216})
+ (6^t + 12^t + 36^t - 40^t + 72^t) S_t(\frac{79}{216}, \frac{475}{1296})\\
&\scriptstyle\;
+ (6^t + 12^t + 36^t - 40^t + 72^t + 216^t + 432^t) S_t(\frac{475}{1296}, \frac{11}{30})
+ 20^t S_t(\frac{11}{30}, \frac{133}{360})
+ (108^t + 216^t) S_t(\frac{259}{648}, \frac{2}{5})\\
&\scriptstyle\;
- (6^t + 18^t + 36^t - 40^t + 108^t + 216^t) S_t(\frac{13}{30}, \frac{281}{648})
- (6^t + 18^t + 36^t - 40^t) S_t(\frac{281}{648}, \frac{313}{720})\\
&\scriptstyle\;
- (6^t + 18^t + 36^t) S_t(\frac{313}{720}, \frac{47}{108})
\mod p,
\end{split}
\end{equation}
with cost $227p/6480<p/28.5$.

\subsection{Heuristic searches}\label{subsecHeuristic}

Our second strategy looks for congruences with small cost, without regard to the number of sums or the complexity of the coefficients in the resulting expression.
This search is designed to investigate \upd{in an empirical way} the question raised by Tanner and Wagstaff \cite{TannerWagstaff} of whether there exist congruences with cost less than $\epsilon p$ for arbitrarily small $\epsilon>0$.
\upd{We provide a general construction that answers this question in the affirmative in Section~\ref{sectionSmallCost}.}

In our heuristic searches, we begin with an initial congruence $\mathcal{C}_0$ and positive integers $D$ and $\lambda$.
We then construct all possible congruences using $\lambda$ subdivision operations, each using a value $d\leq D$.
We select a congruence $\mathcal{C}_1$ of minimal cost from this set, and recursively invoke the algorithm using this as our starting congruence.
This continues for a number of rounds.
We keep $\lambda$ small, often starting with $\lambda=2$, then switching to $\lambda=1$ after a particular number of steps.
We typically select $D=6$ or $D=8$.
This method thus makes locally optimal choices to construct congruences with very low cost.
\upd{Figure~\ref{figGreedyProgress} in Section~\ref{sectionSmallCost} depicts the progress of this method with $\lambda=1$ and $D=8$, beginning with \eqref{eqnB2}.}

Beginning with Vandiver's congruence \eqref{eqnVandiver}, we use this method to construct a congruence for Bernoulli numbers mod $p$ having cost
\[
\frac{5970989p}{241920000} < \frac{p}{40.5}.
\]
The expression we constructed contains $546$ sums, but there is no reason to believe this is optimal.

\section{Searching for Wolstenholme primes}\label{sectionComputationsW}

Using Glaisher's congruence \eqref{eqnGlaisher}, we have that $p\geq5$ is a Wolstenholme prime if and only if $p\mid B_{p-3}$.
We employ the congruences from the prior section to compute the residue of $B_{p-3}$ mod $p$ for primes $p\geq7$, using $k=(p-3)/2$.
In this case, $C_k(a,b,c) \equiv (c^3-a^3-b^3+1)/6$ mod $p$, so for example \eqref{eqnB2} becomes
\[
\textstyle14 B_{p-3} \equiv (1+2^{p-4})S_{p-4}(\frac{1}{6}, \frac{1}{5}) - 2^{p-4}S_{p-4}(\frac{3}{10}, \frac{1}{3}),
\]
that is,
\begin{equation}\label{eqnBB2}
\textstyle112 B_{p-3} \equiv 9S_{p-4}(\frac{1}{6}, \frac{1}{5}) - S_{p-4}(\frac{3}{10}, \frac{1}{3}).
\end{equation}
Here and throughout this section, all congruences are understood to be taken mod $p$.
In the same way, the congruences \eqref{eqnSV}, \eqref{eqnB6}, \eqref{eqnB9}, \eqref{eqnB16}, \eqref{eqnB22}, and \eqref{eqnB30} produce respectively (writing $S(x,y)$ for $S_{p-4}(x,y)$ throughout for economy of space)
\begin{gather}
\textstyle 21 B_{p-3} \equiv S(\frac{1}{6}, \frac{1}{4}),\label{eqnBB1}\\[8pt]
\begin{split}\label{eqnBB6}
458752 B_{p-3} \equiv\ 
&\textstyle -512 S(\frac{3}{20}, \frac{13}{80}) - 520 S(\frac{13}{80}, \frac{1}{6}) + 
 36929 S(\frac{1}{6}, \frac{27}{160})\\
&\textstyle + 36928 S(\frac{27}{160}, \frac{7}{40}) + 
 36864 S(\frac{7}{40}, \frac{1}{5}) - S(\frac{53}{160}, \frac{1}{3}),
 \end{split}\\[8pt]
\begin{split}\label{eqnBB9}
336000 B_{p-3} \equiv\ 
&\textstyle 1000 S(\frac{1}{18}, \frac{11}{180}) + 976 S(\frac{11}{180}, \frac{1}{15}) - 
 101 S(\frac{2}{15}, \frac{5}{36})\\
&\textstyle + 125 S(\frac{7}{36}, \frac{1}{5}) - 24 S(\frac{47}{180}, \frac{4}{15}) - 
 3125 S(\frac{3}{10}, \frac{11}{36})\\
&\textstyle + 24 S(\frac{1}{3}, \frac{61}{180}) + 
 125 S(\frac{13}{36}, \frac{11}{30}) - 24 S(\frac{83}{180}, \frac{7}{15}),
 \end{split}\\[8pt]
 \begin{split}\label{eqnBB16}
72576000 &B_{p-3} \equiv\ 
\textstyle 216000 S(\frac{1}{18}, \frac{11}{180}) + 210816 S(\frac{11}{180}, \frac{7}{108})\\
&\textstyle + 211816 S(\frac{7}{108}, \frac{47}{720}) + 211735 S(\frac{47}{720}, \frac{1}{15})\\
&\textstyle - 25000 S(\frac{1}{10}, \frac{11}{108}) - 21816 S(\frac{2}{15}, \frac{5}{36}) - 648 S(\frac{83}{360}, \frac{25}{108})\\
&\textstyle + 24352 S(\frac{25}{108}, \frac{7}{30}) - 352 S(\frac{4}{15}, \frac{29}{108}) + 648 S(\frac{29}{108}, \frac{97}{360})\\
&\textstyle + 5184 S(\frac{1}{3}, \frac{61}{180}) + 27000 S(\frac{13}{36}, \frac{11}{30}) + 648 S(\frac{11}{30}, \frac{133}{360})\\
&\textstyle + 1000 S(\frac{43}{108}, \frac{2}{5}) - 24919 S(\frac{13}{30}, \frac{313}{720}) - 25000 S(\frac{313}{720}, \frac{47}{108}),
\end{split}\\[8pt]
\begin{split}\label{eqnBB22}
&72576000 B_{p-3} \equiv\ 
\textstyle
216024 S(\frac{1}{18}, \frac{61}{1080})
+ 216000 S(\frac{61}{1080}, \frac{11}{180})\\
&\textstyle
+ 210816 S(\frac{11}{180}, \frac{7}{108})
+ 211816 S(\frac{7}{108}, \frac{47}{720})
+ 211735 S(\frac{47}{720}, \frac{1}{15})\\
&\textstyle
-25000 S(\frac{1}{10}, \frac{11}{108})
-21860 S(\frac{2}{15}, \frac{29}{216})
-21735 S(\frac{29}{216}, \frac{97}{720})
-21816 S(\frac{97}{720}, \frac{5}{36})\\
&\textstyle
-192 S(\frac{119}{540}, \frac{479}{2160})
-195 S(\frac{479}{2160}, \frac{2}{9})
-648 S(\frac{83}{360}, \frac{25}{108})
+ 24352 S(\frac{25}{108}, \frac{7}{30})\\
&\textstyle
+ 3 S(\frac{5}{18}, \frac{601}{2160})
+ 27000 S(\frac{13}{36}, \frac{263}{720})
+ 26919 S(\frac{263}{720}, \frac{79}{216})
+ 27044 S(\frac{79}{216}, \frac{11}{30})\\
&\textstyle
+ 648 S(\frac{11}{30}, \frac{133}{360})
+ 1000 S(\frac{43}{108}, \frac{2}{5})
-24919 S(\frac{13}{30}, \frac{313}{720})
-25000 S(\frac{313}{720}, \frac{47}{108})\\
&\textstyle
+ 192 S(\frac{4}{9}, \frac{241}{540}),
\end{split}
\end{gather}
and
\begin{equation}\label{eqnBB30}
\begin{split}
\scriptstyle
15&\scriptstyle676416000 B_{p-3} \equiv 
46661184 S(\frac{1}{18}, \frac{61}{1080})
+ 46656000 S(\frac{61}{1080}, \frac{11}{180})
+ 45536256 S(\frac{11}{180}, \frac{7}{108})\\
&\scriptstyle
+ 45752256 S(\frac{7}{108}, \frac{47}{720})
+ 45734760 S(\frac{47}{720}, \frac{43}{648})
+ 45735760 S(\frac{43}{648}, \frac{1}{15})
- 5401000 S(\frac{1}{10}, \frac{65}{648})\\
&\scriptstyle
- 5400000 S(\frac{65}{648}, \frac{11}{108})
- 4721885 S(\frac{2}{15}, \frac{173}{1296})
- 4721760 S(\frac{173}{1296}, \frac{29}{216})
- 4694760 S(\frac{29}{216}, \frac{97}{720})\\
&\scriptstyle
- 4712256 S(\frac{97}{720}, \frac{299}{2160})
- 4712904 S(\frac{299}{2160}, \frac{5}{36})
- 41472 S(\frac{119}{540}, \frac{479}{2160})
- 42120 S(\frac{479}{2160}, \frac{2}{9})\\
&\scriptstyle
+ 5184 S(\frac{2}{9}, \frac{241}{1080})
- 139968 S(\frac{83}{360}, \frac{25}{108})
+ 5260032 S(\frac{25}{108}, \frac{151}{648})
+ 5261032 S(\frac{151}{648}, \frac{7}{30})\\
&\scriptstyle
+ 648 S(\frac{5}{18}, \frac{601}{2160})
+ 5832648 S(\frac{13}{36}, \frac{781}{2160})
+ 5832000 S(\frac{781}{2160}, \frac{263}{720})
+ 5814504 S(\frac{263}{720}, \frac{79}{216})\\
&\scriptstyle
+ 5841504 S(\frac{79}{216}, \frac{475}{1296})
+ 5841629 S(\frac{475}{1296}, \frac{11}{30})
+ 139968 S(\frac{11}{30}, \frac{133}{360})
+ 1000 S(\frac{259}{648}, \frac{2}{5})\\
&\scriptstyle
- 5383504 S(\frac{13}{30}, \frac{281}{648})
- 5382504 S(\frac{281}{648}, \frac{313}{720})
- 5400000 S(\frac{313}{720}, \frac{47}{108}).
\end{split}
\end{equation}

It remains to describe how to calculate these sums efficiently.
For this, we first note that
\[
S_{p-4}(x,y) = \sum_{xp<s<yp} s^{p-4} \equiv \sum_{xp<s<yp} s^{-3} \mod p.
\]
We employ the method described in \cite{McIntoshRoettger} (aided there by P.~Montgomery) to compute the latter value.
Corresponding to a sum $S_{p-4}(x,y)$ with $p$ fixed, we define the polynomial
\[
f_{p,x,y}(z) = \prod_{xp<s<yp} (z+s^3) = \sum_{k\geq0} c_k z^k,
\]
and by expansion it follows that $S_{p-4}(x,y) \equiv c_1 c_0^{-1}$ mod $p$.
We may therefore compute the value of $S_{p-4}(x,y)$ mod $p$ with the following simple procedure:
\begin{equation}
\begin{split}\label{eqnComputeS}
&c_0 \leftarrow 1,\; c_1 \leftarrow 0\\
&\textrm{For each integer $s\in(xp,yp)$:}\\
&\quad c_1 \leftarrow c_1 s^3 + c_0 \mod p\\
&\quad c_0 \leftarrow c_0 s^3 \mod p
\end{split}
\end{equation}

We use Nvidia Tesla V100 GPUs to compute these sums.
Each GPU consists of $80$ streaming multiprocessors (SMs), each of which can execute up to $2048$ threads concurrently, although in our computation this needs to reduce to $1024$ threads per SM, owing to our need to work with $64$-bit arithmetic.
Each GPU works in concert with an attached CPU\@.
The CPU is tasked with a range of integers to check.
\upd{We describe our method first for primes $p<6\cdot10^{10}$, and then some improvements implemented for larger primes.}

\subsection{\upd{Method for primes $p<6\cdot10^{10}$}}\label{subsecWSmallPrimes}

\upd{Over this range, computations for one prime are completed before those for another prime are begun, and the GPU is employed to compute each sum in the current congruence in a highly parallel way.}
Given $p$, $\ceiling{xp}$, and $\floor{yp}$, the GPU divides the integers in $(xp, yp)$ into $80\cdot1024$ arithmetic progressions of approximately equal size, one progression for each thread.
Each thread $\theta$ uses the strategy of \eqref{eqnComputeS}, replacing the interval $(xp,yp)$ with an appropriate arithmetic progression $\mathcal{P}_\theta$, to compute the constant coefficient $c_0$ and the linear coefficient $c_1$ mod $p$ for the polynomial
\begin{equation}\label{eqnPolySlice}
\prod_{s\in\mathcal{P}_\theta} (z+s^3).
\end{equation}
The GPU then reduces these $81920$ pairs of values $(c_0,c_1)$ down to $80$ pairs with ten further rounds of parallel calculations.
Each round halves the number of current pairs by combining two pairs $(c_0, c_1)$ and $(c_0', c_1')$ via
\begin{equation}\label{eqnCombinePairs}
(c_1z+c_0)(c_1'z+c_0') \equiv (c_0 c_1' + c_1 c_0') z + c_0 c_0' \mod z^2,
\end{equation}
computing the coefficients mod $p$.
This process ends after ten rounds because only the threads operating in the same block, that is, executing within the same SM, can share memory easily.
After this, the CPU handles the final reduction, using \eqref{eqnCombinePairs} iteratively to combine the remaining $80$ pairs into a single pair, which represents the values of $c_0$ and $c_1$ for the entire sum mod $p$.
Once all the sums in the congruence have been computed, the CPU then combines their values, incorporates the appropriate coefficients, and computes the appropriate inverse mod $p$, to determine the residue of $B_{p-3}$.

It remains to select the congruence employed at each stage of the computation.
While more complicated congruences like \eqref{eqnBB9} and \eqref{eqnBB16} have smaller cost than simpler congruences like \eqref{eqnBB2} and \eqref{eqnBB1}, the congruences with a larger number of sums require greater computational overhead in invoking the GPU multiple times and combining the results of these calls, so the more complex methods are only economical for larger values of $p$.
For example, we found that the nine-term congruence \eqref{eqnBB9} did not offer any savings over the six-term congruence \eqref{eqnBB6} for primes near $1.6\cdot10^{10}$, even though the nominal cost of the former congruence is slightly better.

Consequently, we employed different congruences over our computation.
We employed the simple congruence \eqref{eqnBB1}, with cost $p/12$, for primes $p<1.6\cdot10^{10}$.
We then switched to the six-term congruence \eqref{eqnBB6}, with cost $p/19.2$, for $p<2\cdot10^{10}$.
After this, we used the nine-term congruence \eqref{eqnBB9}, with cost $p/20$, for $p<2.8\cdot10^{10}$, then the $16$-term expression \eqref{eqnBB16}, with cost $p/24$, up to $3.6\cdot10^{10}$, followed by the $22$-term congruence \eqref{eqnBB22}, with cost $3p/80$ up to $5\cdot10^{10}$, and finally the $30$-term congruence \eqref{eqnBB30} with cost $227p/6480<p/28.5$ for $p>5\cdot10^{10}$.
This strategy was not optimal: certainly we would have saved time by employing the two-term congruence \eqref{eqnBB2}, which features cost $p/15$, at an appropriate stage.

Employing GPUs in this way allowed for much faster calculations compared to prior searches.
The congruence \eqref{eqnBB1} was employed in prior searches for Wolstenholme primes \cites{McIntosh,McIntoshRoettger,McIntoshFQ},
and in \cite{McIntoshRoettger} it was reported that determining whether $B_{p-3}\equiv0$ mod $p$ using this congruence for a single prime near $10^9$ required about $4.3$ seconds of CPU time, using MIPS R12000 processors.
Using the same congruence, our GPU implementation requires approximately $3.15$ milliseconds per prime near $10^9$ to compute $B_{p-3}$ mod $p$ on a CPU-GPU pair.
The contrast remains stark when we compare our CPU-GPU implementation with a CPU-only implementation that employs the more modern Intel Xeon Cascade Lake processors: these require about $1.98$ seconds per prime in that range to compute $B_{p-3}$ mod $p$.
We remark also that the GPU method with the slightly more complicated congruence \eqref{eqnBB2} computes $B_{p-3}$ mod $p$ for primes near $10^9$ in an average of just $2.73$ milliseconds per prime, so this two-term congruence is advantageous compared to \eqref{eqnBB1}, even for $p$ near $10^9$.

We employed native arithmetic whenever possible.
For $p<2^{32}$, each multiplication operation mod $p$ produces an intermediate result that fits in a $64$-bit long word, so each such operation required one multiplication instruction and one modular reduction operation. Larger primes required more care. In the general case, in order to multiply $a$ and $b$ mod $p$, with $0\leq a,b<p$, we write $a = a_0 + a_1 2^{32}$ and $b = b_0 + b_1 2^{32}$, so that
\begin{equation}\label{eqnMultiply}
ab = a_0 b_0 + a_1 b_0 2^{32} + b_1(a_0 2^{32} + a_1 2^{64}).
\end{equation}
This allows us to compute this product without arithmetic overflow.
We note that for a given prime $p$ we compute $2^{64}$ mod $p$ just once, using one modular reduction and one addition, so that the amortized cost of computing the expression \eqref{eqnMultiply} mod $p$ is four multiplications, four modular reductions, two bit shifts, and three additions, plus two bitwise ``and'' operations and two additional bit shifts used to create the operands from $a$ and $b$.

We took advantage of certain optimizations that are available over particular ranges.
For example, when computing a sum $S_{p-4}(x,y)$ with $y\leq1/4$ when $p<2^{34}$, then for $xp<s<yp$ the quantity $s^2$ mod $p$ can be calculated natively.
For $p>2^{34}$ we typically needed to use \eqref{eqnMultiply} in its full generality.

For a prime $p$ near $6\cdot10^{10}$, \upd{this strategy} computes the value of $B_{p-3}$ mod $p$ using the $30$-term congruence \eqref{eqnBB30} in approximately $253$ milliseconds.

\subsection{\upd{Method for primes $p>6\cdot10^{10}$}}\label{subsecWLargePrimes}

\upd{%
Three improvements to our strategy implemented for primes larger than $6\cdot10^{10}$ produced significantly faster run times.
\begin{enumerate}
\item Each thread on the GPU needs to compute the $c_0$ and $c_1$ values for \eqref{eqnPolySlice}, and so must compute the values of $s^3$ mod $p$ across the values of an arithmetic progression.
Rather than compute each $s^3$ value independently, however, we can use the value of $s^3$ mod $p$ to produce the value of $(s+\theta)^3$ by adding their difference.
This requires tracking $s^2$, but that value can be maintained with a similar strategy, and likewise for linear terms.
This significantly reduces the number of modular multiplications performed on the GPU, which were fairly expensive operations, implemented using several atomic steps as in \eqref{eqnMultiply}.
This cut our run time by $44\%$ for $p$ near $6\cdot10^{10}$.
\item When computing $a$ mod $p$ on the GPU, we had used an integer remainder operator.
However, instead we can use a floating-point division to determine the approximate quotient, truncate the result, and then use one multiplication and one subtraction to compute the integer remainder.
This is safe provided $p<2^{53}$, and we remain well short of this bound.
This reduced our run times by an additional $89\%$ beyond the improvement from the first item.
\item Rather than treat the primes one by one, splitting the work on each prime across the thousands of threads on the GPU, we can assign each thread on the GPU its own prime, and process $81920$ primes at once.
This reduces the overhead per prime and simplifies our code further (since we now have arithmetic progressions with difference $\theta=1$), but comes with the disadvantage that threads handling larger primes will take longer than those on smaller primes, and we need to wait until all threads complete before assigning a new task to the GPU\@.
On balance, this is a win, gaining an additional $3\%$ in our amortized run time per prime.
\end{enumerate}
Overall, our run times dropped $94\%$ by using these improvements: the amortized cost of testing one prime near $6\cdot10^{10}$ with this strategy using the $30$-term congruence \eqref{eqnBB30} was about $15.3$ milliseconds.
We used this strategy to complete our check to $10^{11}$.
At the end of our run, the time required per prime was about $25.3$ msec.
We can also report that a prime near $10^9$ can be tested with an amortized cost of approximately $0.3$ msec using these improvements.
}

\subsection{\upd{Results}}\label{subsecWResults}
Our \upd{searches} found that no Wolstenholme primes $p$ with $2124679<p<\upd{10^{11}}$ exist.
We record here some ``near misses'': let $\symres{r}{p}$ denote the unique representative of $r$ mod $p$ lying in $(-p/2,p/2]$.
Table~\ref{tableNearMissB50} records the primes $p$ with $10^9<p<\upd{10^{11}}$ where $\abs{\symres{B_{p-3}}{p}} < 50$; a Wolstenholme prime $p$ would of course have $\symres{B_{p-3}}{p} = 0$.
All computations were performed on Gadi, a high-performance computer managed by NCI Australia.
\upd{Beyond $6\cdot10^{10}$, the package \texttt{primesieve} \cite{Walisch} was employed to enumerate primes over intervals.}

\begin{table}[tb]
\caption{Primes $p\in(10^9, \upd{10^{11}})$ for which $\abs{\symres{B_{p-3}}{p}} < 50$.}\label{tableNearMissB50}
\begin{tabular}{c@{\qquad}c}
\begin{tabular}[t]{|r|r|}\hline
\multicolumn{1}{|c|}{$p$} & \multicolumn{1}{c|}{\TS\BS$\symres{B_{p-3}}{p}$}\\\hline
$1025793739$ & $-9$\\
$1029113299$ & $-7$\\
$1939582759$ & $-19$\\
$2139716869$ & $2$\\
$3803691517$ & $13$\\
$8208762073$ & $24$\\
$9267199079$ & $-22$\\
$13581221947$ & $40$\\
$14211360143$ & $-41$\\
$15744104053$ & $-2$\\
\hline
\end{tabular} &
\begin{tabular}[t]{|r|r|}\hline
\multicolumn{1}{|c|}{$p$} & \multicolumn{1}{c|}{\TS\BS$\symres{B_{p-3}}{p}$}\\\hline
$16425136499$ & $7$\\
$21861395221$ & $-11$\\
$22855335949$ & $33$\\
$23345427659$ & $-27$\\
$23543635009$ & $-21$\\
$27827984099$ & $34$\\
$40306537633$ & $42$\\
$44718258259$ & $-6$\\
$56604583391$ & $-25$\\
\upd{$77559232657$} & \upd{$19$}\\
\hline
\end{tabular}
\end{tabular}
\end{table}

\begin{figure}[tb]
\caption{Histogram for $\symres{B_{p-3}}{p}/p$ for $p<\upd{10^{11}}$.}\label{figBernHist}
\begin{center}
\includegraphics[width=3in]{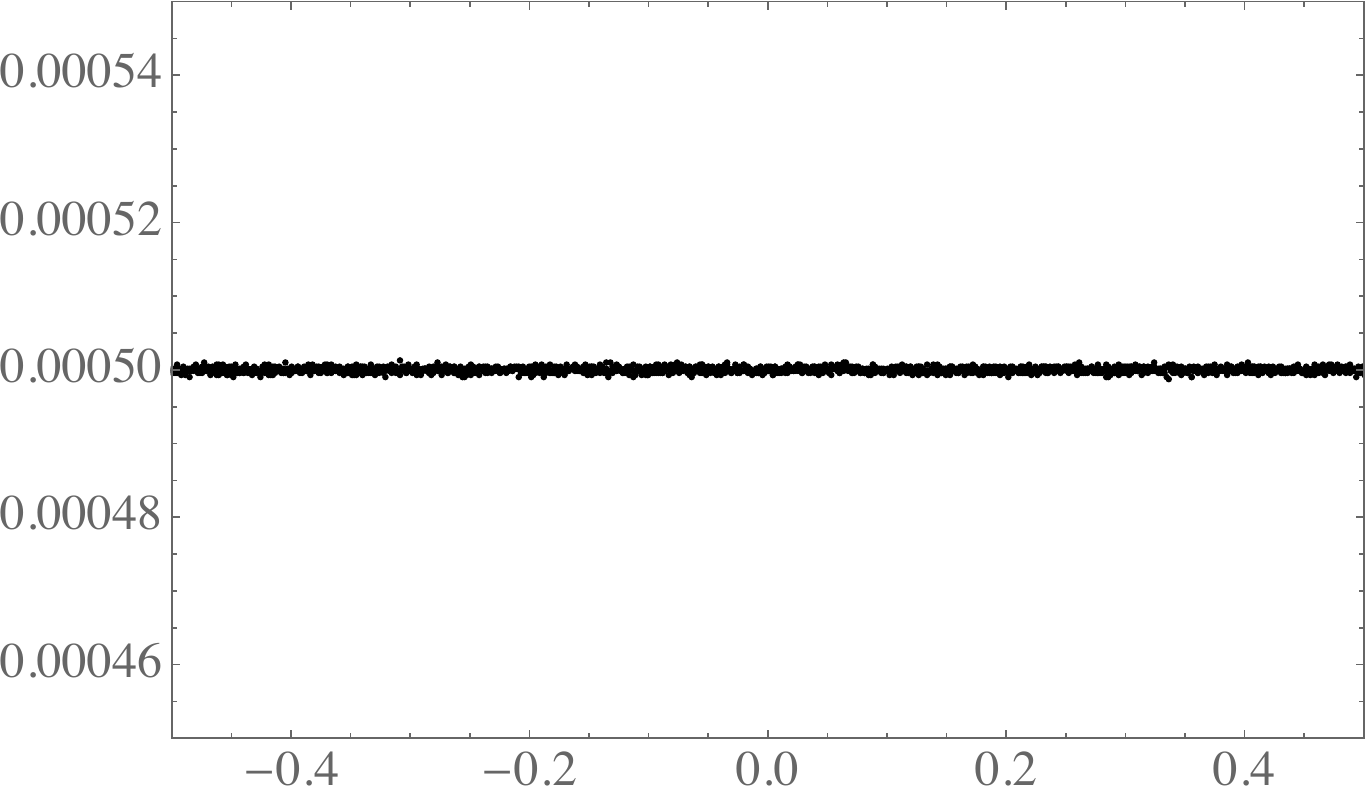}
\end{center}
\end{figure}

We discern no particular bias among the values of $B_{p-3}$ mod $p$.
Figure~\ref{figBernHist} displays a histogram with $2000$ bins of equal size for the values of $\symres{B_{p-3}}{p}/p$ over the \upd{$4118054811$ primes $p$ with $5\leq p<10^{11}$}; a perfectly uniform distribution would display as a horizontal line at $0.0005$.

\section{Congruences for Euler numbers}\label{sectionEuler}

In 1902, Glaisher \cite[\S34]{Glaisher3} established a formula for the residue of Euler numbers modulo a prime $p$.
Using our $S_\ell(x,y)$ notation, this is
\begin{equation}\label{eqnE1}
\textstyle
(-1)^{\frac{p-1}{2}-k} 4^{2k-1} E_{p-1-2k} \equiv S_{p-1-2k}(0,\frac{1}{4})  \mod p
\end{equation}
when $p\geq5$, $2k\leq p-3$ is a positive even integer, and $p-1 \nmid 2k$.
Special cases of this congruence also appear in \cite[p.\ 359]{ELehmer} and \cite[Cor.\ 3.7]{ZHSun08}.
This congruence has cost $p/4$, and was employed in prior searches for Vandiver primes \cite{McIntoshFQ}.
We apply the exhaustive search strategy from Section~\ref{subsecExhaustive}, using similar parameters, to determine some congruences equivalent to \eqref{eqnE1} with smaller cost.

We describe the derivation of four favorable congruences having $m=3$, $5$, $9$, and $16$ terms respectively, which stem from \eqref{eqnE1}, and which were employed in our calculations.
We assume $p\geq5$ and $p-1 \nmid 2k$ throughout.
For our congruence with $m=3$, apply subdivision with $d=2$ to $S_t(0,1/4)$ in \eqref{eqnE1}, and then the same transformation to the resulting term $S_t(0,1/8)$.
After simplifying using reflection and separation as required, we obtain
\begin{equation}\label{eqnE3}
\textstyle
(-1)^k 4^{2k-1} E_t \equiv 4^t S_t(0, \frac{1}{16}) + 2^t S_t(\frac{3}{8}, \frac{7}{16}) + (2^t + 4^t) S_t(\frac{7}{16}, \frac{1}{2}) \mod p,
\end{equation}
where we write $t$ for $p-1-2k$ for brevity.
This congruence has cost $3p/16$.

For $m=5$, we apply two more subdivision transformations on \eqref{eqnE3} with $d=2$, on $S_t(0,1/16)$ and then on $S_t(0,1/32)$, to produce our congruence with cost $9p/64$:
\begin{equation}\label{eqnE5}
\begin{split}
(-1)^k &4^{2k-1} E_t \equiv
\textstyle
16^t S_t(0, \frac{1}{64})
+ 2^t S_t(\frac{3}{8}, \frac{7}{16})
+ (2^t + 4^t) S_t(\frac{7}{16}, \frac{15}{32})\\
&\textstyle
+ (2^t + 4^t + 8^t) S_t(\frac{15}{32}, \frac{31}{64})
+ (2^t + 4^t + 8^t + 16^t) S_t(\frac{31}{64}, \frac{1}{2}) \mod p.
\end{split}
\end{equation}

Next, for $m=9$, we subdivide the $S_t(0,1/128)$ term from \eqref{eqnE5} with $d=2$ and $S_t(3/8,7/16)$ using $d=3$ to produce
\begin{equation}\label{eqnE9}
\begin{split}
(-1)^k &4^{2k-1} E_t \equiv
\textstyle
32^t S_t(0, \frac{1}{128})
+ 6^t S_t(\frac{1}{8}, \frac{7}{48})
+ 6^t S_t(\frac{3}{16}, \frac{5}{24})\\
&\textstyle
+ (2^t + 4^t) S_t(\frac{7}{16}, \frac{11}{24})
+ (2^t + 4^t + 6^t) S_t(\frac{11}{24}, \frac{15}{32})\\
&\textstyle
+ (2^t + 4^t + 6^t + 8^t) S_t(\frac{15}{32}, \frac{23}{48})
+ (2^t + 4^t + 8^t) S_t(\frac{23}{48}, \frac{31}{64})\\
&\textstyle
+ (2^t + 4^t + 8^t + 16^t) S_t(\frac{31}{64}, \frac{63}{128})\\
&\textstyle
+ (2^t + 4^t + 8^t + 16^t + 32^t) S_t(\frac{63}{128}, \frac{1}{2})
\mod p,
\end{split}
\end{equation}
which has cost $43p/384=p/8.93\ldots$\,, so less than half that of \eqref{eqnE1}.
For $m=16$, we subdivide the $S_t(3/16, 5/24)$ and $S_t(7/16, 11/24)$ terms from \eqref{eqnE9} with $d=3$, and $S_t(0,1/128)$, $S_t(1/8, 7/48)$, and $S_t(11/24, 15/32)$ with $d=2$ to yield
\begin{equation}\label{eqnE16}
\begin{split}
&(-1)^k 4^{2k-1} E_t \equiv
\textstyle
64^t S_t(0, \frac{1}{256})
+ (12^t + 18^t) S_t(\frac{1}{16}, \frac{5}{72})
+ 12^t S_t(\frac{5}{72}, \frac{7}{96})\\
&\textstyle\;
+ (6^t + 12^t) S_t(\frac{7}{48}, \frac{11}{72})
+ (6^t + 12^t) S_t(\frac{13}{72}, \frac{3}{16})\\
&\textstyle\;
+ (4^t + 8^t + 12^t) S_t(\frac{11}{48}, \frac{15}{64})
+ 18^t S_t(\frac{19}{72}, \frac{17}{64})\\
&\textstyle\;
+ (4^t + 8^t + 12^t + 18^t) S_t(\frac{17}{64}, \frac{13}{48})
+ 18^t S_t(\frac{19}{48}, \frac{29}{72})
+ 12^t S_t(\frac{41}{96}, \frac{7}{16})\\
&\textstyle\;
+ (2^t + 4^t + 6^t+ 8^t) S_t(\frac{15}{32}, \frac{23}{48})
+ (2^t + 4^t + 6^t + 8^t + 12^t) S_t(\frac{23}{48}, \frac{31}{64})\\
&\textstyle\;
+ (2^t + 4^t + 6^t + 8^t + 12^t + 16^t) S_t(\frac{31}{64}, \frac{35}{72})\\
&\textstyle\;
+ (2^t + 4^t + 8^t + 16^t) S_t(\frac{35}{72}, \frac{63}{128})
+ (2^t + 4^t + 8^t + 16^t + 32^t) S_t(\frac{63}{128}, \frac{127}{256})\\
&\textstyle\;
+ (2^t + 4^t + 8^t + 16^t + 32^t + 64^t) S_t(\frac{127}{256}, \frac{1}{2})
\mod p,
\end{split}
\end{equation}
for a cost of $205p/2304=p/11.23\ldots$\,.
Last, for our relation with $m=24$, we subdivide the terms $S_t(0, \frac{1}{256})$, $S_t(\frac{1}{16}, \frac{5}{72})$, and $S_t(\frac{41}{96}, \frac{7}{16})$ in \eqref{eqnE16} using $d=2$, $4$, and $3$ respectively, and then subdivide $S_t(\frac{41}{288}, \frac{7}{48})$ and $S_t(\frac{15}{32}, \frac{137}{288})$ in the resulting expression with $d=2$ to obtain
\begin{equation}\label{eqnE24}
\begin{split}
&\scriptstyle
(-1)^k 4^{2k-1} E_t \equiv
128^t S_t(0, \frac{1}{512})
+ (48^t+72^t) S_t(\frac{1}{64}, \frac{5}{288})
+ 12^t S_t(\frac{5}{72}, \frac{41}{576})
+ (12^t+72^t) S_t(\frac{41}{576}, \frac{7}{96})\\
&\scriptstyle
+ (6^t+12^t) S_t(\frac{7}{48}, \frac{11}{72})
+ (6^t+12^t) S_t(\frac{13}{72}, \frac{3}{16})
+ 36^t S_t(\frac{3}{16}, \frac{55}{288})
+ (4^t+8^t+12^t) S_t(\frac{11}{48}, \frac{67}{288})\\
&\scriptstyle
+ (4^t+8^t+12^t+48^t+72^t) S_t(\frac{67}{288}, \frac{15}{64})
+ (4^t+8^t+12^t+16^t) S_t(\frac{15}{64}, \frac{137}{576})\\
&\scriptstyle
+ (4^t+8^t+12^t+16^t) S_t(\frac{151}{576}, \frac{19}{72})
+ (4^t+8^t+12^t+16^t+18^t) S_t(\frac{19}{72}, \frac{17}{64})\\
&\scriptstyle
+ (4^t+8^t+12^t+18^t+48^t+72^t) S_t(\frac{17}{64}, \frac{77}{288})
+ (4^t+8^t+12^t+18^t) S_t(\frac{77}{288}, \frac{13}{48})
+ 18^t S_t(\frac{19}{48}, \frac{29}{72})\\
&\scriptstyle
+ 72^t S_t(\frac{41}{96}, \frac{247}{576})
+ (2^t+4^t+6^t+8^t+36^t) S_t(\frac{137}{288}, \frac{23}{48})
+ (2^t+4^t+6^t+8^t+12^t) S_t(\frac{23}{48}, \frac{139}{288})\\
&\scriptstyle
+ (2^t+4^t+6^t+8^t+12^t+48^t+72^t) S_t(\frac{139}{288}, \frac{31}{64})
+ (2^t+4^t+6^t+8^t+12^t+16^t) S_t(\frac{31}{64}, \frac{35}{72})\\
&\scriptstyle
+ (2^t+4^t+8^t+16^t) S_t(\frac{35}{72}, \frac{63}{128})
+ (2^t+4^t+8^t+16^t+32^t) S_t(\frac{63}{128}, \frac{127}{256})\\
&\scriptstyle
+ (2^t+4^t+8^t+16^t+32^t+64^t) S_t(\frac{127}{256}, \frac{255}{512})
+ (2^t+4^t+8^t+16^t+32^t+64^t+128^t) S_t(\frac{255}{512}, \frac{1}{2})
\mod p,
\end{split}
\end{equation}
which has cost $115p/1536=p/13.35\ldots$\,.

We also find that we can derive some better congruences by using another starting relation.
In \cite{McIntoshFQ}, McIntosh proved that
\begin{equation}\label{eqnMacE}
\textstyle
(-1)^k 4^{k-1}(9^k+1) E_{p-1-2k} \equiv \widetilde{S}_{p-1-2k}(0,\frac{1}{6}),
\end{equation}
provided $p\geq5$ and $p \nmid (9^k+1)$, where
\[
\widetilde{S}_\ell(x,y) = \sum_{xp<s<yp} (-1)^s s^\ell \mod p.
\]
It is straightforward to verify that the separation and reflection properties of Proposition~\ref{propTW} hold for $\widetilde{S}_\ell(x,y)$ as well as $S_\ell(x,y)$, but the alternating sign in $\widetilde{S}_\ell(x,y)$ requires an alteration to the subdivision property: if $d$ is a positive integer and $p\nmid d$, then
\[
\widetilde{S}_\ell(x,y) = \begin{cases}
d^\ell \sum_{i=0}^d (-1)^i \widetilde{S}_\ell(\frac{x+i}{d},\frac{y+i}{d}), & \textrm{if $d$ is odd},\\
d^\ell \sum_{i=0}^d (-1)^i S_\ell(\frac{x+i}{d},\frac{y+i}{d}), & \textrm{if $d$ is even}.\\
\end{cases}
\]
We apply a subdivision operation with $d=2$ to the $\widetilde{S}_\ell(x,y)$ expression in \eqref{eqnMacE}, producing (writing $t$ for $p-2k-1$)
\begin{equation}\label{eqnEMac2}
\textstyle
(-1)^k 4^{k-1}(9^k+1) E_t \equiv 2^t\left(S_t(0,\frac{1}{12}) - S_t(\frac{5}{12},\frac{1}{2})\right).
\end{equation}
After this we may apply the transformations of Proposition~\ref{propTW} in the usual way.
Using this strategy, we obtain congruences with costs lower than those we could construct using \eqref{eqnE1} as our base.
We record just one such congruence here, as it is employed in our searches.
Starting with \eqref{eqnEMac2} and applying the sixteen subdivision steps
\begin{equation*}
\begin{split}
&\textstyle
d=2 \textrm{\ on\ }
(0,\frac{1}{2^k\cdot 3}),\; 2\leq k\leq8,\\
&\textstyle
d=3 \textrm{\ on\ }
(\frac{5}{12}, \frac{11}{24}),\;
(\frac{11}{24}, \frac{17}{36}),\;
(\frac{17}{36}, \frac{23}{48}),\;
(\frac{23}{48}, \frac{35}{72}),\\
&\textstyle
d=2 \textrm{\ on\ }
(\frac{5}{36},\frac{11}{72}),\;
(\frac{5}{72},\frac{11}{144}),\\
&\textstyle
d=3 \textrm{\ on\ }
(\frac{61}{144}, \frac{31}{72}),\;
(\frac{133}{288}, \frac{67}{144}),\;
(\frac{205}{432}, \frac{103}{216})
\end{split}
\end{equation*}
in sequence, and simplifying appropriately using reflection and separation after each step, we obtain the following congruence for Euler numbers with cost $27p/512<p/18.96$, which holds for $p\geq5$ and $p\nmid (9^k+1)$:
\begin{equation}\label{eqnE33}
\begin{split}
&\scriptstyle
(-1)^k 4^{k-1}(9^k+1) E_t \equiv
256^t S_t(0,\frac{1}{1536})
-24^t S_t(\frac{5}{144},\frac{11}{288})
-36^t S_t(\frac{61}{432},\frac{31}{216})
-(6^t-12^t) S_t(\frac{11}{72},\frac{133}{864})\\
&\scriptstyle
-(6^t-12^t+72^t) S_t(\frac{133}{864},\frac{67}{432})
-(6^t-12^t) S_t(\frac{67}{432},\frac{17}{108})
-(6^t-12^t+18^t) S_t(\frac{17}{108},\frac{205}{1296})\\
&\scriptstyle
-(6^t-12^t+18^t+108^t) S_t(\frac{205}{1296},\frac{103}{648})
-(6^t-12^t+18^t) S_t(\frac{103}{648},\frac{23}{144})
-(6^t-12^t+18^t-24^t) S_t(\frac{23}{144},\frac{35}{216})\\
&\scriptstyle
-(6^t-12^t+18^t-24^t) S_t(\frac{37}{216},\frac{25}{144})
-(6^t-12^t+18^t) S_t(\frac{25}{144},\frac{113}{648})
-(6^t-12^t+18^t+108^t) S_t(\frac{113}{648},\frac{227}{1296})\\
&\scriptstyle
-(6^t-12^t+18^t) S_t(\frac{227}{1296},\frac{19}{108})
-(6^t-12^t) S_t(\frac{19}{108},\frac{77}{432})
-(6^t-12^t+72^t) S_t(\frac{77}{432},\frac{155}{864})\\
&\scriptstyle
-(6^t-12^t) S_t(\frac{155}{864},\frac{13}{72})
-6^t S_t(\frac{13}{72},\frac{41}{216})
-(6^t+36^t) S_t(\frac{41}{216},\frac{83}{432})
-6^t S_t(\frac{83}{432},\frac{7}{36})\\
&\scriptstyle
-(2^t-4^t+6^t-8^t-12^t) S_t(\frac{35}{72},\frac{421}{864})
-(2^t-4^t+6^t-8^t-12^t+72^t) S_t(\frac{421}{864},\frac{211}{432})\\
&\scriptstyle
-(2^t-4^t+6^t-8^t-12^t) S_t(\frac{211}{432},\frac{47}{96})
-(2^t-4^t+6^t-8^t-16^t-12^t) S_t(\frac{47}{96},\frac{53}{108})\\
&\scriptstyle
-(2^t-4^t+6^t-8^t-16^t-12^t+18^t) S_t(\frac{53}{108},\frac{637}{1296})
-(2^t-4^t+6^t-8^t-16^t-12^t+18^t+108^t) S_t(\frac{637}{1296},\frac{319}{648})\\
&\scriptstyle
-(2^t-4^t+6^t-8^t-12^t-16^t+18^t) S_t(\frac{319}{648},\frac{71}{144})
-(2^t-4^t+6^t-8^t-16^t-12^t+18^t-24^t) S_t(\frac{71}{144},\frac{95}{192})\\
&\scriptstyle
-(2^t-4^t+6^t-8^t-12^t-16^t+18^t-24^t-32^t) S_t(\frac{95}{192},\frac{107}{216})
-(2^t-4^t-8^t-16^t-32^t) S_t(\frac{107}{216},\frac{191}{384})\\
&\scriptstyle
-(2^t-4^t-8^t-16^t-32^t-64^t) S_t(\frac{191}{384},\frac{383}{768})
-(2^t-4^t-8^t-16^t-32^t-64^t-128^t) S_t(\frac{383}{768},\frac{767}{1536})\\
&\scriptstyle
-(2^t-4^t-8^t-16^t-32^t-64^t-128^t-256^t) S_t(\frac{767}{1536},\frac{1}{2})
\mod p.
\end{split}
\end{equation}

Finally, we remark that using the heuristic optimization method of Section~\ref{subsecHeuristic}, with only modest effort we can construct more complicated congruences with cost less than $p/30.4$.
We produced one such congruence having $438$ terms with a greedy search beginning at \eqref{eqnEMac2} and using $\lambda=1$ and $D=8$, iterating over $100$ rounds.
A similar search beginning with \eqref{eqnE1} produced a congruence with $416$ terms and cost near $p/23.6$.

\section{Searching for Vandiver primes}\label{sectionComputationsV}

In order to search for Vandiver primes, we require the case $k=1$ in the congruences of the prior section.
Since $S_{p-3}(x,y)\equiv S_{-2}(x,y)$ mod $p$, we use $t=-2$ in \eqref{eqnE3}, \eqref{eqnE5}, \eqref{eqnE9}, \eqref{eqnE16}, \eqref{eqnE24}, and \eqref{eqnE33} to obtain the following congruences.
Here, we omit the subscript on the $S_t(x,y)$ terms for brevity, and naturally each of these congruences is understood to be taken modulo $p$.
We assume $p\geq7$ throughout.
\begin{gather}
\label{eqnEE3}
-64 E_{p-3} \equiv
\textstyle
S(0, \frac{1}{16}) + 4 S(\frac{3}{8}, \frac{7}{16}) + 5 S(\frac{7}{16}, \frac{1}{2}),\\[8pt]
\begin{split}\label{eqnEE5}
-1024 E_{p-3} \equiv {}
&\textstyle
S(0, \frac{1}{64})
+ 64 S(\frac{3}{8}, \frac{7}{16})
+ 80 S(\frac{7}{16}, \frac{15}{32})\\
&\textstyle
+ 84 S(\frac{15}{32}, \frac{31}{64})
+ 85S(\frac{31}{64}, \frac{1}{2}),
\end{split}\\[8pt]
\begin{split}\label{eqnEE9}
-36864 &E_{p-3} \equiv {}
\textstyle
9 S(0, \frac{1}{128})
+ 256 S(\frac{1}{8}, \frac{7}{48})
+ 256 S(\frac{3}{16}, \frac{5}{24})\\
&\textstyle
+ 2880 S(\frac{7}{16}, \frac{11}{24})
+ 3136 S(\frac{11}{24}, \frac{15}{32})
+ 3280 S(\frac{15}{32}, \frac{23}{48})\\
&\textstyle
+ 3024 S(\frac{23}{48}, \frac{31}{64})
+ 3060 S(\frac{31}{64}, \frac{63}{128})
+ 3069 S(\frac{63}{128}, \frac{1}{2}),
\end{split}\\[8pt]
\begin{split}\label{eqnEE16}
-&1327104 E_{p-3} \equiv {}
\textstyle
81 S(0, \frac{1}{256})
+ 3328 S(\frac{1}{16}, \frac{5}{72})
+ 2304 S(\frac{5}{72}, \frac{7}{96})\\
&\textstyle
+ 11520 S(\frac{7}{48}, \frac{11}{72})
+ 11520 S(\frac{13}{72}, \frac{3}{16})
+ 28224 S(\frac{11}{48}, \frac{15}{64})
+ 1024 S(\frac{19}{72}, \frac{17}{64})\\
&\textstyle
+ 29248 S(\frac{17}{64}, \frac{13}{48})
+ 1024S(\frac{19}{48}, \frac{29}{72})
+ 2304 S(\frac{41}{96}, \frac{7}{16})
+ 118080 S(\frac{15}{32}, \frac{23}{48})\\
&\textstyle
+ 120384 S(\frac{23}{48}, \frac{31}{64})
+ 121680 S(\frac{31}{64}, \frac{35}{72})
+ 110160 S(\frac{35}{72}, \frac{63}{128})\\
&\textstyle
+ 110484 S(\frac{63}{128}, \frac{127}{256})
+ 110565 S(\frac{127}{256}, \frac{1}{2}),
\end{split}\\[8pt]
\begin{split}\label{eqnEE24}
-&5308416 E_{p-3} \equiv {}
\textstyle
81 S(0, \frac{1}{512})
+ 832 S(\frac{1}{64}, \frac{5}{288})
+ 9216 S(\frac{5}{72}, \frac{41}{576})\\
&\textstyle
+ 9472 S(\frac{41}{576}, \frac{7}{96})
+ 46080 S(\frac{7}{48}, \frac{11}{72})
+ 46080 S(\frac{13}{72}, \frac{3}{16})
+ 1024 S(\frac{3}{16}, \frac{55}{288})\\
&\textstyle
+ 112896 S(\frac{11}{48}, \frac{67}{288})
+ 113728 S(\frac{67}{288}, \frac{15}{64})
+ 118080 S(\frac{15}{64}, \frac{137}{576})\\
&\textstyle
+ 118080 S(\frac{151}{576}, \frac{19}{72})
+ 122176 S(\frac{19}{72}, \frac{17}{64})
+ 117824 S(\frac{17}{64}, \frac{77}{288})\\
&\textstyle
+ 116992 S(\frac{77}{288}, \frac{13}{48})
+ 4096 S(\frac{19}{48}, \frac{29}{72})
+ 256 S(\frac{41}{96}, \frac{247}{576})
+ 473344 S(\frac{137}{288}, \frac{23}{48})\\
&\textstyle
+ 481536 S(\frac{23}{48}, \frac{139}{288})
+ 482368 S(\frac{139}{288}, \frac{31}{64})
+ 486720 S(\frac{31}{64}, \frac{35}{72})\\
&\textstyle
+ 440640 S(\frac{35}{72}, \frac{63}{128})
+ 441936 S(\frac{63}{128}, \frac{127}{256})
+ 442260 S(\frac{127}{256}, \frac{255}{512})\\
&\textstyle
+ 442341 S(\frac{255}{512}, \frac{1}{2}),
\end{split}\\[8pt]
 \begin{split}\label{eqnEE33}
&\scriptstyle
477757440 E_{p-3} \equiv {}
-729 S(0,\frac{1}{1536})
+82944 S(\frac{5}{144},\frac{11}{288})
+36864 S(\frac{61}{432},\frac{31}{216})
+995328 S(\frac{11}{72},\frac{133}{864})\\
&\scriptstyle
+1004544 S(\frac{133}{864},\frac{67}{432})
+995328 S(\frac{67}{432},\frac{17}{108})
+1142784 S(\frac{17}{108},\frac{205}{1296})
+1146880 S(\frac{205}{1296},\frac{103}{648})\\
&\scriptstyle
+1142784 S(\frac{103}{648},\frac{23}{144})
+1059840 S(\frac{23}{144},\frac{35}{216})
+1059840 S(\frac{37}{216},\frac{25}{144})
+1142784 S(\frac{25}{144},\frac{113}{648})\\
&\scriptstyle
+1146880 S(\frac{113}{648},\frac{227}{1296})
+1142784 S(\frac{227}{1296},\frac{19}{108})
+995328 S(\frac{19}{108},\frac{77}{432})
+1004544 S(\frac{77}{432},\frac{155}{864})\\
&\scriptstyle
+995328 S(\frac{155}{864},\frac{13}{72})
+1327104 S(\frac{13}{72},\frac{41}{216})
+1363968 S(\frac{41}{216},\frac{83}{432})
+1327104 S(\frac{83}{432},\frac{7}{36})\\
&\scriptstyle
+9206784 S(\frac{35}{72},\frac{421}{864})
+9216000 S(\frac{421}{864},\frac{211}{432})
+9206784 S(\frac{211}{432},\frac{47}{96})
+9020160 S(\frac{47}{96},\frac{53}{108})\\
&\scriptstyle
+9167616 S(\frac{53}{108},\frac{637}{1296})
+9171712 S(\frac{637}{1296},\frac{319}{648})
+9167616 S(\frac{319}{648},\frac{71}{144})
+9084672 S(\frac{71}{144},\frac{95}{192})\\
&\scriptstyle
+9038016 S(\frac{95}{192},\frac{107}{216})
+7978176 S(\frac{107}{216},\frac{191}{384})
+7966512 S(\frac{191}{384},\frac{383}{768})
+7963596 S(\frac{383}{768},\frac{767}{1536})\\
&\scriptstyle
+7962867 S(\frac{767}{1536},\frac{1}{2}).
\end{split}
\end{gather}

We may compute each sum $S_{p-3}(x,y)$ in the same manner as \upd{our Wolstenholme searches}, employing the same polynomial method by computing the linear and constant coefficients of the polynomial $\prod_{xp<s<yp} (z+s^2)$, similar to the procedure of \eqref{eqnComputeS}.
These sums were again computed using GPUs\upd{, using a strategy similar to that of Section~\ref{subsecWSmallPrimes} for primes $p<4\cdot10^{10}$, and utilizing the improvements of Section~\ref{subsecWLargePrimes} for larger primes.
The first improvement, tracking differences across an arithmetic progression when computing powers, yielded a more modest $28\%$ gain in this case, since we require squares here instead of cubes, but our overall gain over all three improvements was quite similar to the Wolstenholme case, at $93\%$.
Near $4\cdot10^{10}$, the amortized cost of computing $E_{p-3}$ mod $p$ for a single prime using \eqref{eqnEE33} dropped from approximately $195$ milliseconds before the improvements to about $14.0$ msec with them.
Near $10^{11}$, the average time per prime measured about $34.1$ msec.}

We employed the congruence \eqref{eqnEE3} for the search range $p<10^9$, followed by \eqref{eqnEE5} for $10^9 < p < 2^{32}$, then \eqref{eqnEE9} for $2^{32}<p<2^{33}$, then \eqref{eqnEE16} for $2^{33}<p<2^{34}$, followed by \eqref{eqnEE24} for $2^{34}<p<2\cdot10^{10}$, and finally \eqref{eqnEE33} for $p>2\cdot10^{10}$.

Our searches found no new Vandiver primes larger than $p=1062232319$.
Table~\ref{tableNearMissE50} records the primes $p$ with $10^9<p<\upd{10^{11}}$ where $\abs{\symres{E_{p-3}}{p}} < 50$; a Vandiver prime $p$ has $\symres{E_{p-3}}{p} = 0$.
Computations were again completed at NCI Australia\upd{, and \texttt{primesieve} \cite{Walisch} employed for enumerating primes above $4\cdot10^{10}$}.

\begin{table}[tb]
\caption{Primes $p\in(10^9, \upd{10^{11}})$ for which $\abs{\symres{E_{p-3}}{p}} < 50$.}\label{tableNearMissE50}
\begin{tabular}{c@{\qquad}c}
\begin{tabular}[t]{|r|r|}\hline
\multicolumn{1}{|c|}{$p$} & \multicolumn{1}{c|}{\TS\BS$\symres{E_{p-3}}{p}$}\\\hline
$1062232319$ & $0$\\
$1348936931$ & $-17$\\
$1352698411$ & $17$\\
$1836806681$ & $-15$\\
$2114780851$ & $-2$\\
$2161739347$ & $32$\\
$2264214119$ & $-38$\\
$2978890751$ & $-35$\\
$3700821251$ & $-23$\\
$10158743171$ & $-49$\\
$10179358499$ & $-12$\\
$14884379297$ & $-40$\\
\hline
\end{tabular} &
\begin{tabular}[t]{|r|r|}\hline
\multicolumn{1}{|c|}{$p$} & \multicolumn{1}{c|}{\TS\BS$\symres{E_{p-3}}{p}$}\\\hline
$17380814081$ & $5$\\
$18044797027$ & $-10$\\
$18642203467$ & $34$\\
$23177794127$ & $-48$\\
$36652898767$ & $-9$\\
$36830964851$ & $6$\\
\upd{$42574679263$} & \upd{$-2$}\\
\upd{$52794781757$} & \upd{$-24$}\\
\upd{$58999725809$} & \upd{$43$}\\
\upd{$66004284323$} & \upd{$19$}\\
\upd{$73607557631$} & \upd{$7$}\\
\hline
\end{tabular}
\end{tabular}
\end{table}

\begin{figure}[tb]
\caption{Histogram for $\symres{E_{p-3}}{p}/p$ for $p<\upd{10^{11}}$.}\label{figEulerHist}
\begin{center}
\includegraphics[width=3in]{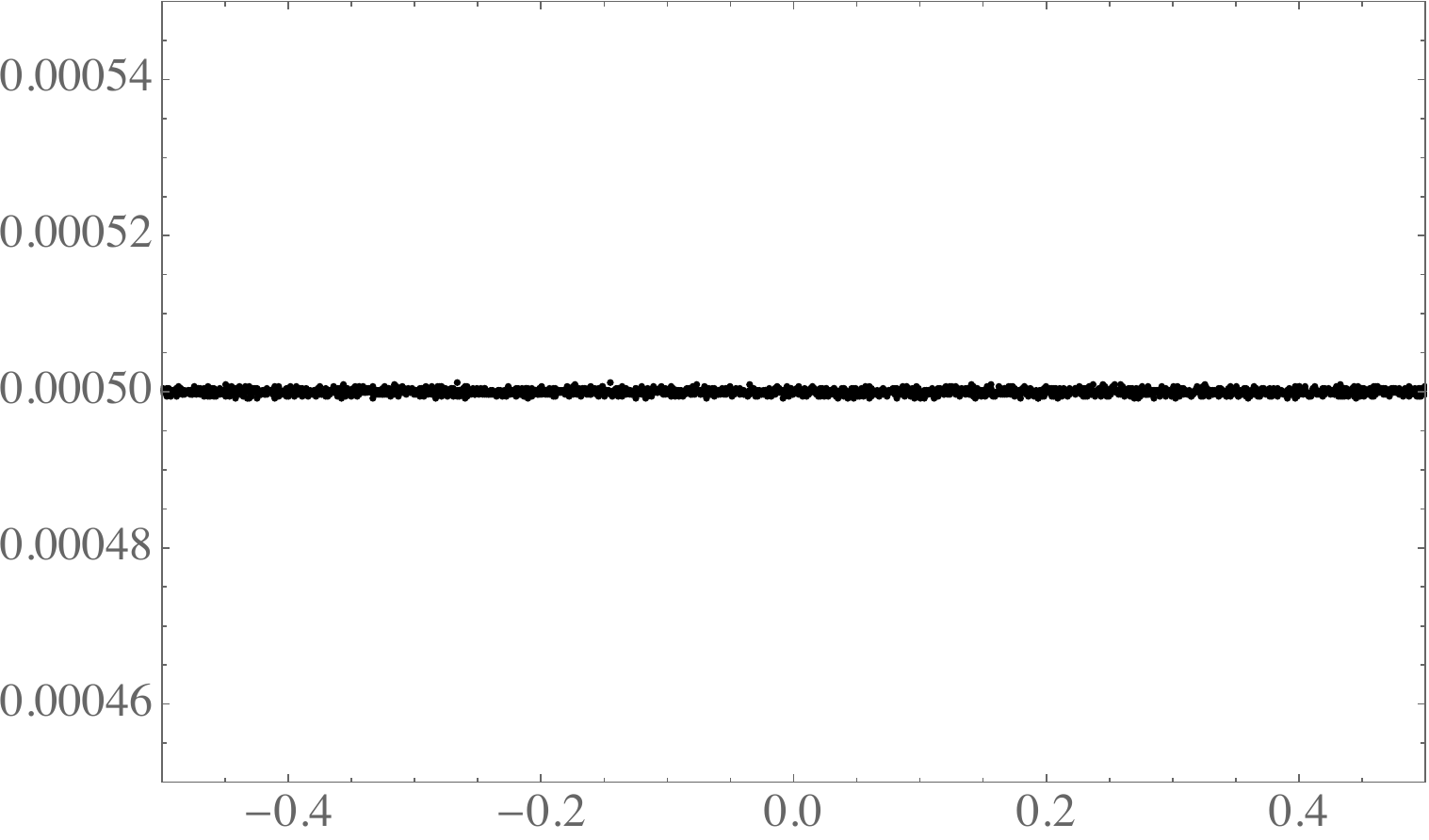}
\end{center}
\end{figure}

As with the Bernoulli numbers, we observe no particular bias among the values of $E_{p-3}$ mod $p$.
Figure~\ref{figEulerHist} displays a histogram with $2000$ equal-sized bins for the values of $\symres{E_{p-3}}{p}/p$ over the primes $p$ with $5\leq p<\upd{10^{11}}$.

\section{\upd{Congruences with arbitrarily small cost}}\label{sectionSmallCost}

\upd{%
In \cite{TannerWagstaff}, Tanner and Wagstaff raised the question of whether there exist congruences for computing Bernoulli numbers involving linear combinations of the sums $S_\ell(x,y)$ having cost less than $\epsilon p$, for any $\epsilon>0$.
One may ask the same question regarding the computation of residues for Euler numbers.
We answer both questions in the affirmative here, by showing that any such sum has an equivalent formulation with arbitrarily small cost.
}

\begin{thm}\label{thmSmallCost}
\upd{Let $x<y$ be rational numbers in $[0,1]$, let $p$ be an odd prime number that does not divide the denominators of $x$ and $y$, and let $\ell$ be an integer.
Given $\epsilon>0$, there exist disjoint intervals $(x_1,y_1)$, \ldots, $(x_m,y_m)$ in $[0,1]$ and integers $c_1$, \ldots, $c_m$, with the $x_i$, $y_i$, and $c_i$ independent of $p$, so that
\[
S_\ell(x,y) \equiv \sum_{i=1}^m c_i S_\ell(x_i,y_i) \mod p
\]
where the cost of the right side is less than $\epsilon p$.}
\end{thm}

\begin{proof}
\upd{%
Let $q$ be a common denominator for $x$ and $y$ with the property that $y-x=n/q$ with $4\mid n$.
Set $A\subseteq\{0,\ldots,q-1\}$ so that
\[
S_\ell(x,y) = \sum_{a\in A} S_\ell\left(\frac{a}{q},\frac{a+1}{q}\right)
\]
by the separation property of Proposition~\ref{propTW}.
Partition $A$ into two subsets $A_1$ and $A_2$ of equal size.
For each $a\in A_1$, use the subdivision property with $d=q$ to write
\begin{equation}\label{eqnSCA}
S_\ell\left(\frac{a}{q},\frac{a+1}{q}\right) \equiv q^\ell \sum_{i=0}^{q-1} S_\ell\left(\frac{a+iq}{q^2}, \frac{a+iq+1}{q^2}\right) \mod p.
\end{equation}
For each $b\in A_2$, use the separation property to write
\begin{equation}\label{eqnSCB}
S_\ell\left(\frac{b}{q},\frac{b+1}{q}\right) = \sum_{j=0}^{q-1} S_\ell\left(\frac{bq+j}{q^2}, \frac{bq+j+1}{q^2}\right).
\end{equation}
For each pair $(a,b)\in A_1\times A_2$, the right sides of \eqref{eqnSCA} and \eqref{eqnSCB} have one term in common, when $i=b$ and $j=a$.
These $\abs{A_1}\abs{A_2}=n^2/4$ terms are disjoint, so the cost of our new expression is
\[
\frac{pn}{q} - \frac{pn^2}{4q^2} = \frac{pn}{q}\left(1-\frac{n}{4q}\right).
\]
The transformed sum involves $n(q-n/4)$ intervals, each of size $1/q^2$, and since $4\mid n(q-n/4)$ we can iterate this process.
Using the fact that $x(1-x/4)$ is increasing on $[0,1]$, a straightforward inductive argument shows that after $k$ iterations the cost of our expression is less than $4p/k$, and we may choose $k$ large enough so that this is $<\epsilon p$.
}
\end{proof}

\upd{By applying this procedure to \eqref{eqnSV} and \eqref{eqnE1}, it follows immediately that there exist congruences for Bernoulli numbers $B_k$ and Euler numbers $E_k$ mod $p$ with arbitrarily small cost, assuming $2k\nmid(p-1)$.}

\begin{figure}[tb]
\caption{\upd{Progress achieved by the method of Section~\ref{subsecHeuristic} for} \upd{Bernoulli congruences (using $\lambda=1$, $D=8$), and a logarithmic fit.}
\upd{Each data point represents a congruence having $m$ terms and cost $p/r$.}}\label{figGreedyProgress}
\begin{center}
\includegraphics[width=3.25in]{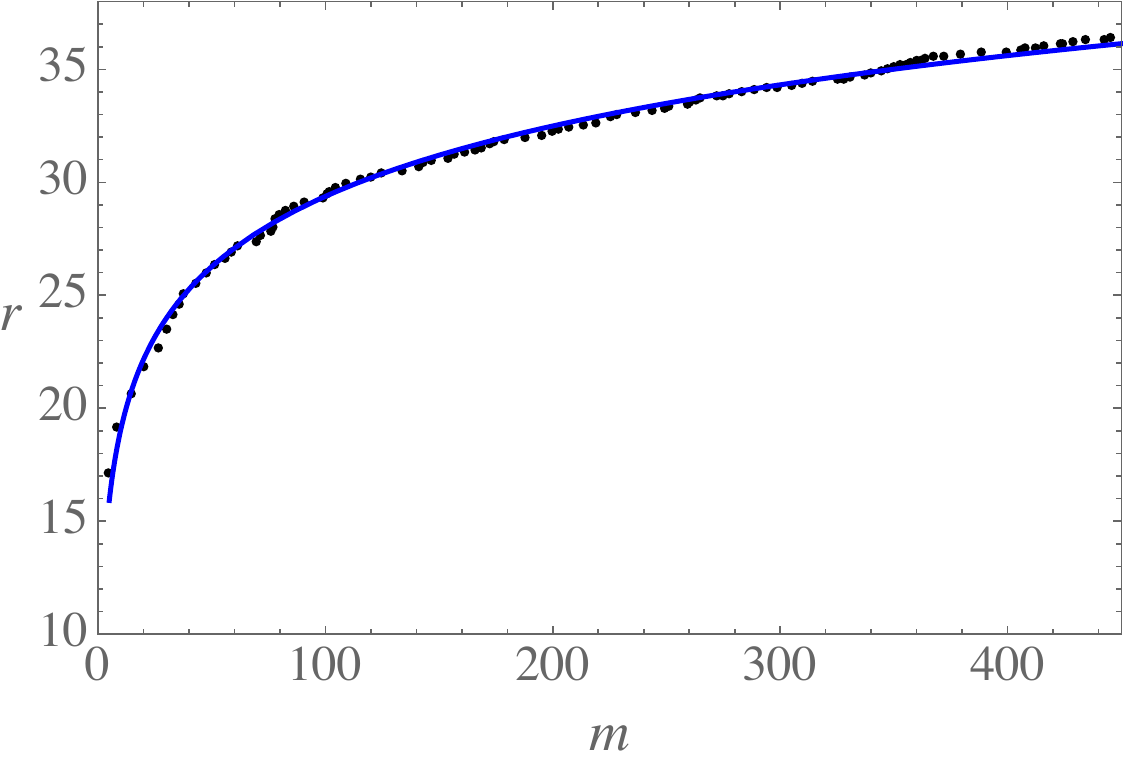}
\end{center}
\end{figure}

\upd{We remark also that the convergence rate obtained by this method is rather slow: since there are $m=\Omega(q^{2^k})$ terms in the sum after $k$ iterations, the cost is $O(p/\log \log m)$.
Perhaps an alternative method could supply an improved convergence rate: empirical evidence from our computations in Section~\ref{subsecHeuristic} and the end of Section~\ref{sectionEuler} suggests that $O(p/\log m)$ is possible, and we conjecture this to be the case.
Figure~\ref{figGreedyProgress}, for example, shows the progress made by iterating this greedy method, beginning with \eqref{eqnB2}, and using parameters $\lambda=1$ and $D=8$.
Here, we plot the point $(m,r)$ if a congruence produced had $m$ terms and cost $p/r$.
The best-fit curve of the form $a+b\log m$ is also displayed: here $a=8.68$ and $b=4.49$.
}

\upd{We provide a few basic bounds related to derivations of congruences in the next section.}

\section{The transformation graph}\label{sectionBounds}

While not directly applicable to the search for Wolstenholme and Vandiver primes, this section is devoted to developing ideas of the Tanner and Wagstaff machinery that generates new congruence relations for Bernoulli numbers and Euler numbers. These congruence relations are crucial to the aforementioned search and any insight into their derivation could be useful in future work.

We define a \textit{derived congruence relation} as one that is generated from another congruence relation of the type (\ref{eqnTW1}) using the transformations of Proposition \ref{propTW}---\textit{separation}, \textit{reflection}, and \textit{subdivision}. To conveniently represent these transformations, we will denote separation by $\sigma_f$, where $0 \le f \le 1$ denotes the fraction by which we split the interval into two, i.e., if $0\leq x<y\leq 1$, then a transformation $\sigma_f$ applied on $S_\ell(x,y)$ will result in 
\begin{equation*}
    S_\ell(x,y) = S_\ell\left(x,x+f\cdot (y-x)\right)+S_\ell \left(x+f\cdot (y-x),y\right).
\end{equation*}
Similarly, reflection is denoted by $\rho$ and subdivision into $d$ terms is denoted by $\tau_d$. We will also use $\iota$ to denote the \textit{identity} transformation which will be used only to indicate termination of a chain of transformations.

\begin{figure}[tb]
    \caption{The transformation graph for \eqref{eqnB9}.}
    \centering
    \includegraphics[width=4.9in]{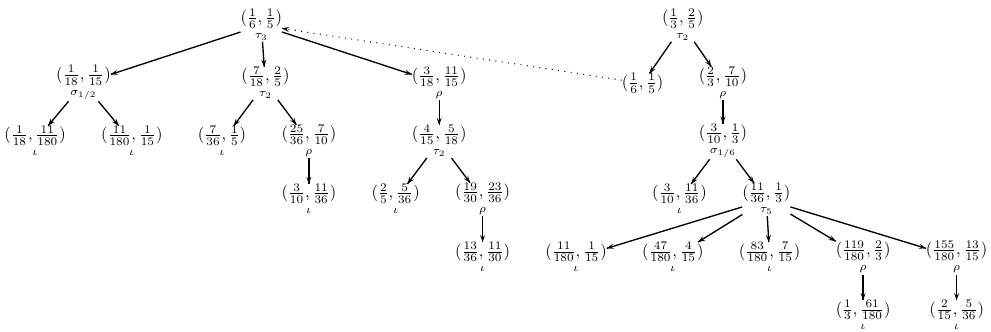}
    \label{m9congruencegraph}
\end{figure}

Starting from (\ref{eqnTW1}), we can draw a tree of transformations corresponding to each of the $\beta:=\floor{b/2}$ terms, where identity/terminal ($\iota$), separation ($\sigma_f$), reflection ($\rho$), and subdivision ($\tau_d$) result in one, two, one, and $d$ branches, respectively. The union of the $\beta$ trees will be referred to as a \textit{transformation graph}. This graph can be traversed in a depth-first-search manner, sequentially from the first tree to the $\beta$th tree, to give a unique string identifying the series of transformations. For example, Figure \ref{m9congruencegraph} depicts the transformations needed to derive (\ref{eqnB9}) from (\ref{eqnTW1}) with $b = 5$ in the form of a graph with two trees corresponding to the two terms in the original congruence relation. The string of transformations corresponding to this is
\begin{equation*}
    \left(\tau_3\sigma_{\frac{1}{2}}\iota\iota\tau_2\iota\rho\iota\rho\tau_2\iota\rho\iota\right) \left(\tau_2 \tau_3\sigma_{\frac{1}{2}}\iota\iota\tau_2\iota\rho\iota\rho\tau_2\iota\rho\iota \rho\sigma_{\frac{1}{6}}\iota\tau_5\iota\iota\iota\rho\iota\rho\iota\right).
\end{equation*}
The final congruence (\ref{eqnB9}) can be obtained by merging the leaf nodes, after multiplying them by the appropriate coefficients.

Using the idea of a transformation graph, we can now prove a basic bound relating to the question posed by Tanner and Wagstaff in \cite{TannerWagstaff} on the existence of congruences with cost less than $\epsilon p$. Assume that the desired cost of a derived congruence is less than or equal to $\epsilon p$ and the least number of terms possible for such a congruence is $T$. Here, a \textit{term} is an expression of the form $\alpha S_\ell(x,y)$ in the derived congruence, $\alpha$ being the appropriate coefficient generated during the transformations. Let the total number of leaf nodes (which are then merged into $T$ terms) in the transformation graph be $T'$. The derived congruence has the form
\begin{equation*}
    C_k(2,b,b+1)B_{2k} \equiv \sum_{i = 1}^T \alpha_i S_{2k-1}(x_i,y_i) \mod p.
\end{equation*}
Since the desired cost is at most $\epsilon p$, we have $\sum_{i = 1}^T (y_i-x_i) \le \epsilon$ and therefore, 
\begin{equation} \label{eqnFbound}
    y_i - x_i \le \epsilon
\end{equation}
for each $i$ and
\begin{equation} \label{eqnFTbound}
    y_j - x_j \le \frac{\epsilon}{T}
\end{equation}
for some $j$.

Trivially, $T' \ge \beta$. The cost of each term in (\ref{eqnTW1}) is $mp/b(b+1)$. If $mp/b(b+1) > \epsilon p$, we need to apply a string of transformations to obtain terms with reduced cost. Among the various transformations, only separation and subdivision yield terms with lower cost than the original. To reduce the cost of each term while generating the least number of extra terms, it is clear that one cannot do better than applying a single subdivision operation, $\tau_{d(m)}$, such that (from (\ref{eqnFbound}))
\begin{equation*}
    \frac{mp}{b(b+1)n(m)} \le \epsilon p,
\end{equation*}
where $d(m)$ is an integer depending on $m$. The value $d(m)$ might be interpreted as the number of leaf nodes in the $m$th tree corresponding to the $m$th term of (\ref{eqnTW1}). Note that other sequences of transformations might exist that generate the same number of terms with the required cost but for our purposes, this assumption is sufficient. 

Therefore, we obtain $d(m) \ge m/\epsilon b(b+1)$. In other words, $\lceil m/\epsilon b(b+1) \rceil$ is a lower bound on the number of branches (and hence, leaf nodes) that need to be generated in each tree in order to satisfy (\ref{eqnFbound}). Therefore, 
\begin{equation} \label{ineqT'1}
    T' \ge \sum_{m=1}^\beta \left \lceil \frac{m}{\epsilon b(b+1)} \right \rceil \ge \max \left\{\beta, \frac{\beta(\beta+1)}{2 \epsilon b(b+1)} \right\}.
\end{equation}
From (\ref{eqnFTbound}), we know that for some $m_0$, we have $d(m_0) \ge m_0T/\epsilon b(b+1)$, or
\begin{equation} \label{ineqT}
    T \le \frac{d(m_0) \epsilon b(b+1)}{m_0} \le N \epsilon b(b+1)
\end{equation}
where $N := \max \limits_m \{d(m)\}$ is the maximum number of leaf nodes in any tree in the graph. In fact, we can refine (\ref{ineqT'1}) slightly:
\begin{equation} \label{ineqT'2}
    T' \ge \max \left\{\beta-1+N, \frac{\beta(\beta-1)}{2 \epsilon b(b+1)} + N \right\}.
\end{equation}
Combining (\ref{ineqT}) and (\ref{ineqT'2}) yields the following result.

\begin{prop} \label{propproportion}
$\displaystyle\frac{T}{T'} \leq \frac{N\epsilon b(b+1)}{\max \left\{\beta-1+N, \frac{\beta(\beta-1)}{2 \epsilon b(b+1)} + N \right\}}$.
\end{prop}
The fraction $T/T'$ is the ratio of distinct leaf nodes to all leaf nodes, so Proposition~\ref{propproportion} implies that, depending on $b$, $N$, and $\epsilon$, a certain proportion of leaf nodes must occur more than once.
In particular, for fixed $b$ and $N$, as $\epsilon\to0^+$  the ratio $T/T'$ approaches $0$.

\section*{Acknowledgements}

We thank Karl Dilcher, Lars Hesselholt, and Richard McIntosh for helpful correspondence.
We also thank NCI Australia and UNSW Canberra for computational resources.
This research was undertaken with the assistance of resources and services from the National Computational Infrastructure (NCI), which is supported by the Australian Government.

\bibliographystyle{amsplain}

\begin{bibdiv}
\begin{biblist}

\bib{AebiCairns}{article}{
   author={Aebi, C.},
   author={Cairns, G.},
   title={Wolstenholme again},
   journal={Elem. Math.},
   volume={70},
   date={2015},
   number={3},
   pages={125--130},
   issn={0013-6018},
   review={\MR{3372069}},
}

\bib{Agoh}{article}{
   author={Agoh, T.},
   title={Congruences related to the Ankeny-Artin-Chowla conjecture},
   journal={Integers},
   volume={16},
   date={2016},
   pages={paper no.\ A12, 30 pp.},
   review={\MR{3475253}},
}

\bib{BCEM93}{article}{
   author={Buhler, J.},
   author={Crandall, R.},
   author={Ernvall, R.},
   author={Mets\"{a}nkyl\"{a}, T.},
   title={Irregular primes and cyclotomic invariants to four million},
   journal={Math. Comp.},
   volume={61},
   date={1993},
   number={203},
   pages={151--153},
   issn={0025-5718},
   review={\MR{1197511}},
}

\bib{Carlitz54}{article}{
   author={Carlitz, L.},
   title={Note on irregular primes},
   journal={Proc. Amer. Math. Soc.},
   volume={5},
   date={1954},
   pages={329--331},
   issn={0002-9939},
   review={\MR{61124}},
}

\bib{CosgraveDilcher}{article}{
   author={Cosgrave, J.~B.},
   author={Dilcher, K.},
   title={On a congruence of Emma Lehmer related to Euler numbers},
   journal={Acta Arith.},
   volume={161},
   date={2013},
   number={1},
   pages={47--67},
   issn={0065-1036},
   review={\MR{3125151}},
}

\bib{Edwards}{book}{
   author={Edwards, H.~M.},
   title={Fermat's Last Theorem: A Genetic Introduction to Algebraic Number Theory},
   series={Grad. Texts in Math.},
   volume={50},
   publisher={Springer, New York},
   date={1996},
   pages={xvi+410},
   isbn={0-387-90230-9},
   isbn={0-387-95002-8},
   review={\MR{1416327}},
}

\bib{Ernvall79}{article}{
   author={Ernvall, R.},
   title={Corrigenda: ``Cyclotomic invariants and $E$-irregular primes''},
   journal={Math. Comp.},
   volume={33},
   date={1979},
   number={145},
   pages={433},
   issn={0025-5718},
   review={\MR{514840}},
}

\bib{EM78}{article}{
   author={Ernvall, R.},
   author={Mets\"{a}nkyl\"{a}, T.},
   title={Cyclotomic invariants and $E$-irregular primes},
   journal={Math. Comp.},
   volume={32},
   date={1978},
   number={142},
   pages={617--629},
   issn={0025-5718},
   review={\MR{482273}},
}

\bib{Gardiner}{article}{
   author={Gardiner, A.},
   title={Four problems on prime power divisibility},
   journal={Amer. Math. Monthly},
   volume={95},
   date={1988},
   number={10},
   pages={926--931},
   issn={0002-9890},
   review={\MR{1581615}},
}

\bib{Glaisher1}{article}{
   author={Glaisher, J.~W.~L.},
   title={Congruences relating to the sums of products of the first $n$ numbers and to other sums of products},
   journal={Quart. J. Pure Appl. Math.},
   volume={31},
   date={1900},
   pages={1--35},
}

\bib{Glaisher2}{article}{
   author={Glaisher, J.~W.~L.},
   title={On the residues of the sums of products of the first $p-1$ numbers, and their powers, to modulus $p^2$ or $p^3$},
   journal={Quart. J. Pure Appl. Math.},
   volume={31},
   date={1900},
   pages={321--353},
}

\bib{Glaisher3}{article}{
   author={Glaisher, J.~W.~L.},
   title={On the residues of the sums of the inverse powers of numbers in arithmetic progression},
   journal={Quart. J. Pure Appl. Math.},
   volume={32},
   date={1901},
   pages={271--305},
}

\bib{Gut}{article}{
   author={Gut, M.},
   title={Euler'sche Zahlen und Klassenanzahl des K\"{o}rpers der 4$\ell$-ten Einheitswurzeln},
   journal={Comment. Math. Helv.},
   volume={25},
   date={1951},
   pages={43--63},
   issn={0010-2571},
   review={\MR{41173}},
}

\bib{HHO}{article}{
   author={Hart, W.},
   author={Harvey, D.},
   author={Ong, W.},
   title={Irregular primes to two billion},
   journal={Math. Comp.},
   volume={86},
   date={2017},
   number={308},
   pages={3031--3049},
   issn={0025-5718},
   review={\MR{3667037}},
}

\bib{Johnson75}{article}{
   author={Johnson, W.},
   title={Irregular primes and cyclotomic invariants},
   journal={Math. Comp.},
   volume={29},
   date={1975},
   pages={113--120},
   issn={0025-5718},
   review={\MR{376606}},
}

\bib{ELehmer}{article}{
   author={Lehmer, E.},
   title={On congruences involving Bernoulli numbers and the quotients of Fermat and Wilson},
   journal={Ann. of Math. (2)},
   volume={39},
   date={1938},
   number={2},
   pages={350--360},
   issn={0003-486X},
   review={\MR{1503412}},
}

\bib{LPMP15}{article}{
   author={Luca, F.},
   author={Pizarro-Madariaga, A.},
   author={Pomerance, C.},
   title={On the counting function of irregular primes},
   journal={Indag. Math. (N.S.)},
   volume={26},
   date={2015},
   number={1},
   pages={147--161},
   issn={0019-3577},
   review={\MR{3281697}},
}

\bib{McIntosh}{article}{
   author={McIntosh, R.~J.},
   title={On the converse of Wolstenholme's theorem},
   journal={Acta Arith.},
   volume={71},
   date={1995},
   number={4},
   pages={381--389},
   issn={0065-1036},
   review={\MR{1339137}},
}

\bib{McIntoshFQ}{article}{
   author={McIntosh, R.~J.},
   title={Congruences involving Euler numbers and power sums},
   journal={Fibonacci Quart.},
   volume={58},
   date={2020},
   number={4},
   pages={328--333},
}

\bib{McIntoshRoettger}{article}{
   author={McIntosh, R.~J.},
   author={Roettger, E.~L.},
   title={A search for Fibonacci-Wieferich and Wolstenholme primes},
   journal={Math. Comp.},
   volume={76},
   date={2007},
   number={260},
   pages={2087--2094},
   issn={0025-5718},
   review={\MR{2336284}},
}

\bib{MestrovicSurvey}{article}{
   author={Me\v{s}trovi\'c, R.},
   title={Wolstenholme's theorem: its generalizations and extensions in the last hundred and fifty years (1862--2012)},
   note={arXiv.1111.3057v2},
   date={25 Dec.\ 2011},
   pages={31 pp.},
}

\bib{Mestrovic14}{article}{
   author={Me\v{s}trovi\'{c}, R.},
   title={A search for primes $p$ such that the Euler number $E_{p-3}$ is divisible by $p$},
   journal={Math. Comp.},
   volume={83},
   date={2014},
   number={290},
   pages={2967--2976},
   issn={0025-5718},
   review={\MR{3246818}},
}

\bib{SelfridgePollack}{article}{
   author={Selfridge, J.~L.},
   author={Pollack, B.~W.},
   title={Fermat's last theorem is true for any exponent up to 25,000},
   journal={Notices Amer. Math. Soc.},
   volume={11},
   number={1},
   date={1964},
   pages={97},
   note={Abstract no. 608-138, Annual Meeting in Miami and Coral Gables, FL, Jan. 23-26, 1964},
}

\bib{StaffordVandiver}{article}{
   author={Stafford, E.~T.},
   author={Vandiver, H.~S.},
   title={Determination of some properly irregular cyclotomic fields},
   journal={Proc. Natl. Acad. Sci. USA},
   volume={16},
   date={1930},
   number={2},
   pages={139--150},
}

\bib{ZHSun08}{article}{
   author={Sun, Z.-H.},
   title={Congruences involving Bernoulli and Euler numbers},
   journal={J. Number Theory},
   volume={128},
   date={2008},
   number={2},
   pages={280--312},
   issn={0022-314X},
   review={\MR{2380322}},
   doi={10.1016/j.jnt.2007.03.003},
}

\bib{ZWSun11}{article}{
   author={Sun, Z.-W.},
   title={Super congruences and Euler numbers},
   journal={Sci. China Math.},
   volume={54},
   date={2011},
   number={12},
   pages={2509--2535},
   issn={1674-7283},
   review={\MR{2861289}},
}

\bib{TannerWagstaff}{article}{
   author={Tanner, J.~W.},
   author={Wagstaff, S.~S., Jr.},
   title={New congruences for the Bernoulli numbers},
   journal={Math. Comp.},
   volume={48},
   date={1987},
   number={177},
   pages={341--350},
   issn={0025-5718},
   review={\MR{866120}},
}

\bib{UspenskyHeaslet}{book}{
   author={Uspensky, J.~V.},
   author={Heaslet, M.~A.},
   title={Elementary Number Theory},
   publisher={McGraw-Hill, New York},
   date={1939},
   pages={x+484},
   review={\MR{0000236}},
}

\bib{VTW}{article}{
   author={Van Der Poorten, A.~J.},
   author={te Riele, H.~J.~J.},
   author={Williams, H.~C.},
   title={Computer verification of the Ankeny-Artin-Chowla conjecture for all primes less than $100\,000\,000\,000$},
   journal={Math. Comp.},
   volume={70},
   date={2001},
   number={235},
   pages={1311--1328},
   review={\MR{1709160} and MR{1933835}},
   note={Corrigenda and addition, ibid. \textbf{72} (2003), no. 241, 521--523},
}

\bib{Vandiver37}{article}{
   author={Vandiver, H.~S.},
   title={On Bernoulli's numbers and Fermat's last theorem},
   journal={Duke Math. J.},
   volume={3},
   date={1937},
   number={4},
   pages={569--584},
   issn={0012-7094},
   review={\MR{1546011}},
}

\bib{Vandiver40}{article}{
   author={Vandiver, H.~S.},
   title={Note on Euler number criteria for the first case of Fermat's last theorem},
   journal={Amer. J. Math.},
   volume={62},
   date={1940},
   pages={79--82},
   issn={0002-9327},
   review={\MR{1231}},
}

\bib{Wagstaff}{article}{
   author={Wagstaff, S.~S., Jr.},
   title={The irregular primes to $125000$},
   journal={Math. Comp.},
   volume={32},
   date={1978},
   number={142},
   pages={583--591},
   issn={0025-5718},
   review={\MR{491465}},
}

\bib{Walisch}{article}{
   author={\upd{W}\upd{alisch}, \upd{K.}},
   title={\upd{Primesieve, ver. 7.6}},
   date={\upd{Jan. 7, 2021}},
   note={\upd{\url{https://github.com/kimwalisch/primesieve}}},
}

\bib{Washington}{book}{
   author={Washington, L.~C.},
   title={Introduction to Cyclotomic Fields},
   series={Grad. Texts in Math.},
   volume={83},
   edition={2},
   publisher={Springer, New York},
   date={1997},
   pages={xiv+487},
   isbn={0-387-94762-0},
   review={\MR{1421575}},
}

\bib{Wolstenholme}{article}{
   author={Wolstenholme, J.},
   title={On certain properties of prime numbers},
   journal={Quart. J. Pure Appl. Math.},
   volume={5},
   date={1862},
   pages={35--39},
   issn={0033-5606},
}

\end{biblist}
\end{bibdiv}

\end{document}